\documentclass[12pt, leqno]{article}
\usepackage{amssymb,amsmath,latexsym}
\evensidemargin =\oddsidemargin

\font\teneufm=eufm10
\font\seveneufm=eufm7
\font\fiveeufm=eufm5
\newfam\eufmfam
\textfont\eufmfam=\teneufm
\scriptfont\eufmfam=\seveneufm
\scriptscriptfont\eufmfam=\fiveeufm

\let\goth\mathfrak

\def\gg{\goth g}
\def\ff{\goth f}
\def\gs{\goth s}

\def\gk{\goth k}
\def\gl{\goth l}

\def\gp{\goth p}

\def\gn{\goth n}

\def\gq{\goth q}
\def\gz{\goth z}
\def\go{\goth o}
\def\ga{\goth a}

\def\pp{\mbox{\bf p}}

\def\ad{\mbox{ad}}

\def\beq{\begin{equation}}
\def\eeq{\end{equation}}
\def\bea{\begin{eqnarray}}
\def\eea{\end{eqnarray}}
\def\beas{\begin{eqnarray*}}
\def\eeas{\end{eqnarray*}}

\def\cplus{\hbox{$\supset${\raise1.05pt\hbox{\kern -0.55em
${\scriptscriptstyle +}$}}\ }}

\DeclareMathOperator{\Span}{Span}
\DeclareMathOperator{\tr}{tr}
\DeclareMathOperator{\Hom}{Hom}
\DeclareMathOperator{\Aut}{Aut}

\DeclareMathOperator{\Der}{Der}
\DeclareMathOperator{\Ctd}{Ctd}

\DeclareMathOperator{\End}{End}
\DeclareMathOperator{\Id}{Id}

\DeclareMathOperator{\Spec}{Spec}

\DeclareMathOperator{\Ad}{Ad}
\DeclareMathOperator{\Spin}{Spin}
\DeclareMathOperator{\Cliff}{Cliff}
\DeclareMathOperator{\Gr}{Gr}

\newtheorem{theorem}[equation]{Theorem}

\newtheorem{lemma}[equation]{Lemma}

\newtheorem{proposition}[equation]{Proposition}

\def\Z{\mathbb Z}
\def\C{\mathbb C}

\title{Automorphisms and twisted loop algebras of finite dimensional 
simple Lie superalgebras}
\author{Dimitar Grantcharov\thanks{Supported by a Pacific Institute 
for the Mathematical Sciences postdoctoral fellowship} \, and  Arturo 
Pianzola\thanks{Supported by NSERC operating grant A9343}}
\date{}
\begin{document}

\maketitle

\begin{abstract}
We describe the structure of the algebraic group of automorphisms of all simple finite dimensional Lie superalgebras. Using this and \'etale cohomology considerations, we list all different isomorphism classes of the corresponding twisted loop superalgebras.
\end{abstract}

2000 MSC: Primary 17B67. Secondary 17B40, 20G15. 

\section{Introduction}

Twisted loop algebras were introduced by V. Kac to provide explicit realizations of all affine Kac-Moody Lie 
algebras (see Ch. 8 of \cite{K1} for details). The basic ingredients of his construction are a finite 
dimensional simple complex Lie algebra $\gg$, and an automorphism 
$\sigma$ of $\gg$ of finite order. The loop algebras then appear as 
complex (infinite dimensional) Lie subalgebras $L(\gg, \sigma) 
\subset \gg \otimes \C [z^{\pm 1}]$. Deciding when two of these loop algebras 
are isomorphic is a delicate question. This can be done using the invariance of the underlying Cartan
matrices (which requires conjugacy of Cartan subalgebras) or cohomological methods (see below). The final outcome
is that $\gg$ is an invariant of $L(\gg, \sigma)$, and that the $\sigma$'s may be 
chosen to come from the different symmetries of the Dynkin diagram 
(up to conjugation).

In \cite{S} the author sets out to recreate this 
construction in the case when $\gg$ is  a finite dimensional 
simple complex Lie superalgebra. Therein one finds some rather intricate
 arguments (crucial to the present work, see (a) below) that 
determine the structure of 
the (abstract) quotient group $\mbox{ Aut}\, \gg /  G_0$, where 
$G_0$ is the group of ``inner automorphisms'' of $\gg$. (The group $G_0$
 is related to the group of 
automorphisms of the even part of $\gg$, but it is not the connected 
component of $  \mbox{Aut}\, \gg$.) The author shows that the elements of $G_{0}$ 
lead to (untwisted) loop algebras, and states that a 
given list of coset representatives of $  \mbox{ Aut}\, \gg /  G_0$ accounts for 
the remaining different isomorphism classes of twisted loop algebras. 
The issue of the invariance of $\gg$ is not raised.\footnote{Serganova has informed us that she  did have case-by-case
 computational proofs of some of these results in her unpublished Ph.D thesis.
In any event, we believe that the cohomological point of view explains the result better and sheds considerably
 more insight  into these questions.}

In \cite{P} and \cite{ABP} it is shown that the study of isomorphism 
classes of loop algebras of arbitrary algebras is a problem of Galois 
cohomology. More precisely, one views $L(\gg, \sigma)$ as an algebra 
over the ring $R = k[z^{\pm m}]$ (where $m$ is the order of $\sigma$ and $k$ is our base field) 
and shows that if  $S = k[z^{\pm 1}]$ then $L(\gg, \sigma)\otimes_R S \simeq (\gg \otimes 
R)\otimes_R S$ as $S$-algebras. Thus $L(\gg, \sigma)$ is an $S/R$ 
form of $\gg \otimes R$. Since the extension $S/R$ is \'etale (in 
fact finite Galois), we can attach to $L(\gg, \sigma)$ an element 
(isomorphism class of a torsor)
$$
[L(\gg, \sigma)] \in 
H_{\acute et}^1(\Spec R, \mbox{\bf Aut}\, \gg _R)
$$
(see 
\cite{P}  for details). The study of this $H^1$ is enough to 
understand the isomorphism classes of loop algebras as {\it algebras 
over $R$}. The passage from $R$ to the base field $k$ is done by 
bringing the centroid of $\gg$ into the picture (\cite{ABP}, \S4). This is the approach that we will follow, 
and as can be seen 
from the above discussion then, there are three key ingredients that 
must be looked at:

\medskip

(a) The nature of the {\it algebraic group}  $\mbox{\bf
Aut}\, \gg $.

(b) The computation of $H_{\acute et}^1(\Spec R,\mbox{\bf 
Aut}\, \gg_R )$.

(c) The supercentroid of $\gg$.

\medskip
(a) is dealt 
with in \S\ref{proof1}. With (b) in mind we describe for each $\gg$ 
(according to the classification due to Kac) an exact sequence of 
algebraic groups 
$$
1 \to \mbox{\bf G}^0 \to \mbox{\bf Aut} \, \gg 
\to \mbox{F}_k \to 1
$$
where $\mbox{\bf G}^0 $ is the connected component of the identity of  $\mbox{\bf Aut} \, \gg$ and 
$\mbox{F}$ is finite. The work of Serganova, which essentially tells us what   $\mbox{\bf Aut} \, \gg(k)$ is,
plays here a crucial role.\footnote{Care must be taken though, to avoid the use of the ``analytic'' part of her arguments using the exponential
 function over $\C$. This is done by introducing suitable Cartan subgroups. }

(b) This is done in \S\ref{proof2}. With the aid of (a) and Grothendieck's work on the algebraic fundamental group, 
we obtain natural correspondences  
$$
H_{\acute et}^1(\Spec R, \mbox{\bf Aut}\, \gg_R ) \simeq H^1(\Spec 
R, \mbox{F}_R ) \simeq \mbox{Conjugacy classes of F}.
$$
This crucial finiteness result is analogous to one of the deep theorems of classical Galois cohomology 
over fields (\cite{CG}, Ch. III, Th\'eor\`eme 4.4). Its validity rests ultimately on the very special nature of the Dedekind ring
$k[t^{\pm 1}]$. The parameterization of $R$-isomorphism 
classes of loop algebras by conjugacy classes of $\mbox{F}$ will follow from this result.

(c) This is the extent of \S\ref{centro}. Once $\gg$ is shown to be central,
we reason as in \cite{ABP}.

\medskip
{\it Acknowledgements:} We would like to thank the referee for suggesting the proof of Proposition \ref{central}
 which is substantially shorter than the original one.

\section{Conventions and notation}

{\it Henceforth $k$ will denote an 
algebraically closed field of characteristic zero.} The category of 
associative commutative unital $k$-algebras will be denoted by 
$k${\it -alg}.  By $\Z_n$ we denote the group $\Z / n\Z$.  We adhere to the convention of writing algebra filtrations $A =  A^l \supset A^{l+1} \supset... $ with upper indexes and algebra gradings $A = \oplus_i A_i$ with lower indexes. 

Let $R \in k${\it -alg}. Recall   that a {\it Lie 
superalgebra over} $R$ is a $\Z_2$-graded  $R$-module $L=L_{\bar{0}} 
\oplus L_{\bar{1}} $ together with a multiplication $[\cdot ,\cdot ]: L \times 
L \to L$ satisfying the following three axioms:

\medskip
(i) $[\cdot , \cdot ]$ is 
$R$-bilinear,

(ii) $[x,y] = - (-1)^{{\deg x}{\deg y}}[y,x]$ ($\deg u 
= i$ if $u \in L_{\bar{i}}$), 

(iii) $[x,[y,z]]+(-1)^{\deg x(\deg y 
+ \deg z)}[y,[z,x]]+ (-1)^{\deg z(\deg x + \deg y)}[z,[x,y]] =0$.

\medskip

An {\it automorphism} of $L$ is an $R$-module automorphism $\varphi$ 
of $L$ satisfying:
 
(i) $\varphi$ stabilizes $L_{\bar{0}}$ and  
$L_{\bar{1}}$ (i.e. $\varphi$ is homogeneous of degree $\bar{0}$),

(ii) 
$\varphi([x,y]) = [\varphi(x), \varphi(y)]$ for all $x,y \in L$.

Endomorphisms are defined in a similar fashion. Note that in some works (see \cite{K2} for example) automorphisms and endomorphisms
are not supposed to preserve the $\Z_2$-grading. For details on all simple finite dimensional Lie superalgebras over $k$ 
we refer the reader to  \cite{K2} and \cite{Sch}. In the present paper we will use the notations of \cite{Pen}.

\smallskip
{\it Assume henceforth that $\gg = \gg_{\bar{0}} \oplus \gg_{\bar{1}}$ is 
a finite dimensional simple Lie superalgebra over $k$.}

\smallskip
By base change 
$\gg(R) = \gg_{\bar{0}}(R) \oplus \gg_{\bar{1}}(R)$ becomes a Lie 
superalgebra over $R$. Consider the $k$-group {\bf Aut}$\, \gg : R 
\to Aut_{R\mbox{-}Ls} \gg(R)$  where the right hand side denotes the group of 
automorphisms of $\gg(R)$ as a Lie superalgebra over $R$. It is easy 
to see that  {\bf Aut}$\, \gg$ is a closed subgroup of {\bf 
GL}$(\gg)$, hence affine (see \S \ref{Gen} below).

\section{Loop algebras} \label{loop}

\setcounter{equation}{0}

We fix once and for all a 
compatible set $\{ \zeta_m \}$ of primitive $m$-th roots of unity in 
$k$. (Thus $\zeta_{lm}^l = \zeta_m$). Let $\sigma $ be an automorphism 
of $\gg$ of finite period $m$. We then get a $\Z_m$-grading
$\oplus_{i = 0}^{m-1}\gg_{\bar{i}}$ of $\gg$ via the 
eigenspace decomposition $\gg_{\bar{i}} = \{ x \in \gg \; | 
\; \sigma(x) = \zeta_m^i x\}$. Define 
$$
L_m(\gg, \sigma):= 
\oplus_{i \in \Z}\gg_{\bar{i}}\otimes z^i \subset \gg \otimes 
k[z^{\pm 1}].
$$
Since by definition $\sigma$ preserves the even and odd parts of 
$\gg$, we see that $L_m(\gg, \sigma)$ is in natural way an (infinite 
dimensional) Lie superalgebra over $k$. An easy calculation shows that up to 
isomorphism $L_m(\gg, \sigma)$ does not depend on the choice of the 
period $m$ of $\sigma$, therefore by a harmless abuse in notation we will 
henceforth simply write $L(\gg, \sigma)$ instead. This is the {\it loop algebra} of 
the pair $(\gg, \sigma)$. (The reason for writing periods and not 
orders is that one can then deal with the isomorphism question 
$L(\gg, \sigma) \simeq L(\gg, \tau)$ in easier terms).

\section{Statements of the main theorems} \label{main}

\begin{theorem} \label{th1}
There exists an exact sequence of linear algebraic $k$-groups
$$
1 \to \mbox{\bf G}^0 \to \mbox{\bf 
Aut}\, \gg \to \mbox{\rm F}_k \to 1.
$$
where $\mbox{\bf G}^0$ (the connected component of the identity of $\mbox{\bf Aut}\, \gg$)  
and $\mbox{\rm F}_k$ (finite and constant) are given by the following table.\footnote{The $\mbox{\bf N}$'s  appearing in the Cartan type series are unipotent groups. Their specific description is given in \S \ref{cartan}}

\end{theorem}

\bigskip
\mbox{Table: Automorphism groups of the 
simple Lie superalgebras.}\footnote{Although the Lie superalgebras $D(1)$ and $\go \gs \gp (4|2)$ are isomorphic we will treat them separately given that the two
 different viewpoints yield different information.}
$$
\begin{tabular}{|c|c|c|c|}
\hline

${\gg}$ & $\mbox{\bf G}^0 = \mbox{\bf Aut}^0\, \gg $ & \mbox{\rm F}$_k$& 
{\bf Split} \\ \hline
$\gs \gl (m|n)$ & $(\mbox{\bf SL}_m \times 
\mbox{\bf SL}_n \times \mbox{\bf G}_m)/ (\boldsymbol{\mu}_m \times \boldsymbol{\mu}_n) $& 
$\Z_{2,k}$ & \mbox{No} \\ \hline

$\gp \gs \gl (n|n), n>2$ & 
$(\mbox{\bf SL}_n \times \mbox{\bf SL}_n \times \mbox{\bf G}_m)/ 
(\boldsymbol{\mu}_n \times \boldsymbol{\mu}_n) $& $\Z_{2,k} \times \Z_{2,k}$ & \mbox{No} \\ 
\hline

$\gp \gs \gl (2|2)$ & $(\mbox{\bf SL}_2 \times \mbox{\bf 
SL}_2 \times \mbox{\bf SL}_2 )/ (\boldsymbol{\mu}_2 \times \boldsymbol{\mu}_2) $& $\Z_{2,k} $ & 
\mbox{Yes} \\ \hline

$\gs \pp (n)$ & $(\mbox{\bf SL}_n \times 
\mbox{\bf G}_m)/ \boldsymbol{\mu}_n $ & $\mbox{1}_k$ & \mbox{Yes} \\ \hline

$\gp 
\gs \gq (n)$ & $\mbox{\bf PGL}_n $ & $ \Z_{4,k}$  & \mbox{Yes} \\ 
\hline

$\go \gs \gp (2l|2n)$ & $(\mbox{\bf SO}_{2l} \times \mbox{\bf 
Sp}_{2n}) / \boldsymbol{\mu}_2$& $\Z_{2,k} $ & \mbox{Yes} \\ \hline

$\go \gs \gp 
(2l+1|2n)$ & $\mbox{\bf SO}_{2l+1} \times \mbox{\bf Sp}_{2n}$ & 
$\mbox{1}_k$   & \mbox{Yes} \\ \hline

$F(4)$ & $(\mbox{\bf Spin}_{7} \times \mbox{\bf SL}_{2})/\boldsymbol{\mu}_2$ & $\mbox{1}_k$   & \mbox{Yes} \\ \hline

$G(3)$ & $ \mbox{\bf G}_{2} \times \mbox{\bf SL}_{2}$ & $\mbox{1}_k$   & \mbox{Yes} \\ \hline

$D(\alpha), \alpha^3 \neq 1, \alpha \neq 
-(2)^{\pm 1}$ &$ (\mbox{\bf SL}_{2} \times \mbox{\bf SL}_{2} \times 
\mbox{\bf SL}_{2}) / (\boldsymbol{\mu}_2 \times \boldsymbol{\mu}_2)$ & $\mbox{1}_k$   & 
\mbox{Yes} \\ \hline

$D(\alpha), \alpha \in \{1, -2, -1/2\}$ & $( 
\mbox{\bf SL}_{2} \times \mbox{\bf SL}_{2} \times \mbox{\bf SL}_{2}) 
/ (\boldsymbol{\mu}_2 \times \boldsymbol{\mu}_2)$ & $ \Z_{2,k} $  & \mbox{Yes} \\ \hline

$D(\alpha), \alpha^3 = 1, \alpha \neq 1 $ &$ (\mbox{\bf SL}_{2} 
\times \mbox{\bf SL}_{2} \times \mbox{\bf SL}_{2}) / (\boldsymbol{\mu}_2 \times 
\boldsymbol{\mu}_2)$ & $\Z_{3,k} $  & \mbox{Yes} \\ \hline

$W(n)$ & $\mbox{\bf 
N}_{W(n)} \rtimes \mbox{\bf GL}_{n} $& $\mbox{1}_k$  & \mbox{Yes} \\ 
\hline

$S(n)$ &$ \mbox{\bf N}_{S(n)} \rtimes \mbox{\bf GL}_{n} $& 
$\mbox{1}_k$  & \mbox{Yes} \\ \hline

$S'(2l)$ &$ \mbox{\bf 
N}_{S'(2l)} \rtimes \mbox{\bf SL}_{2l} $& $\mbox{1}_k$  & \mbox{Yes} 
\\ \hline

$H(2l)$ &$ \mbox{\bf N}_{H(2l)} \rtimes ((\mbox{\bf 
SO}_{2l} \times \mbox{\bf G}_m)/\boldsymbol{\mu}_{2})$ & $\Z_{2,k}$  & \mbox{Yes} 
\\ \hline

$H(2l+1)$ &$ \mbox{\bf N}_{H(2l+1)} \rtimes (\mbox{\bf 
SO}_{2l+1} \times \mbox{\bf G}_m)$ &  $\mbox{1}_k$   & \mbox{Yes} \\ 
\hline
\end{tabular}
$$

\medskip
\noindent
{\bf Note.} It is natural to consider the supergroup (superspace) of automorphisms of $\gg$ in the spirit of Bernstein, Leites, Manin, etc. (see \cite{DM} and \cite{Man}). More precisely, one defines a functor $ \mbox{\bf Aut}^s \, \gg : R \to Aut_{R\mbox{-}Ls} \gg (R)$ from the category $k${\it -salg} of associative unital $k$-superalgebras to the category of sets. We can then relate the purely even superspace $(\mbox{\bf Aut}^s \, \gg)_{\rm rd} $ to the (purely even) affine algebraic group  $\mbox{\bf Aut} \, \gg $ above (recall that for a superspace $(M, {\cal O}_M)$, by $M_{\rm rd}$ we denote $(M, {\cal O}_M/{\cal J}_M)$, where ${\cal J}_M = {\cal O}_{M,1} +  {\cal O}_{M,1}^2$). The explicit 
 structure of $\mbox{\bf Aut}^s \, \gg $ is not needed in the present paper and will be left for future work.

The Lie algebra of the algebraic group $ \mbox{\bf Aut} \, \gg$ is comprised of the (parity preserving) derivations of the algebra $\gg$. These are precisely the
 even derivations of $\gg$ considered as a Lie superalgebra. In this way  we recover from Theorem \ref{th1} the list of the even derivations $(\Der \gg)_0$ of all simple finite-dimensional Lie superalgebras $\gg$
(for the full list see Appendix A in \cite{Pen}).

\medskip
\begin{theorem} \label{cohomo}
Let $R = k[t^{\pm 1}]$. The canonical map $H_{\acute et}^1(R, \mbox{\bf Aut}\, \gg_R) 
\to H_{\acute et}^1(R,\mbox{\rm F}_R )$ is bijective.
\end{theorem}

\begin{theorem} \label{th3}
For $i = 1,2$ let $\gg_i$ be a finite dimensional 
simple Lie superalgebra over $k$, and $\sigma_i \in  \mbox{\bf Aut}\, \gg_i$ an automorphism of finite order.  If $L(\gg_1, \sigma_1) \simeq L(\gg_2, \sigma_2)$ then $\gg_1 \simeq \gg_2$ (all isomorphisms as $k$-Lie superalgebras).

\end{theorem}

\begin{theorem} \label{th3.5}
The following two conditions are equivalent:

(i) $L(\gg, \sigma_1)$ and  $L(\gg, \sigma_2)$ are isomorphic as $k$-Lie superalgebras,

(ii) $\bar{\sigma}_1$ is conjugate to either $\bar{\sigma}_2$ or 
$\bar{\sigma}^{-1}_2$ in $ \mbox{\rm F}:=  \mbox{\rm F}_k(k)$, where  $\; \bar{}:  
\mbox{\bf Aut}\, \gg(k) \to  \mbox{\rm F}$ is the surjective 
(abstract) group isomorphism arising from Theorem \ref{th1}.
\end{theorem}

\begin{theorem} \label{th4}
 Let $R = k[t^{\pm 1}]$ and consider the $R$-Lie 
superalgebras  $\gg \otimes_k R$. Let $\ff$ be an $R$-form 
of $\gg \otimes_k R$; namely  $\ff$ is an $R$-Lie superalgebra such that  $\ff 
\otimes_R S \simeq (\gg \otimes_k R)\otimes_R S \simeq \gg \otimes_k 
S$ as Lie superalgebras for some faithfully flat and finitely 
presented (commutative associative unital) $R$-algebra $S$. Then $\ff 
\simeq L(\gg, \sigma)$ for some $\sigma$ and conversely. Furthermore, 
$S$ can be taken to be finite \'etale (in fact Galois).

\end{theorem}

\section{Generalities about Lie superalgebras and 
algebraic groups} \label{Gen}

\setcounter{equation}{0}
 
Because of our present needs, we will use the 
``functorial approach'' to algebraic groups. This facilitates 
immensely the explicit definition of morphisms and quotients. We 
briefly recall how this goes (see \cite{DG}, \cite{W}, and specially  
Chapter 6 of \cite{KMRT}). 

A $k${\it -group} is a functor from $k${\it -alg} into the category of 
groups. These form a category where morphisms are natural 
transformations. Thus a morphism $f: \mbox{\bf G} \to \mbox{\bf H}$ of 
$k$-groups is nothing but a functorial collection of abstract group homomorphisms 
$f_R : \mbox{\bf G}(R) \to \mbox{\bf H}(R)$, for $R \in k${\it -alg}. 
In particular if $\varphi : R \to R'$ is a $k$-algebra homomorphism then the 
diagram 
$$
\begin{array}{ccc}
\mbox{\bf G}(R)           &  \stackrel{f_R}{\longrightarrow}   &    
\mbox{\bf H}(R)     \\
\mbox{\bf G}(\varphi) \downarrow   \hspace{.5in}  
&                                    &   \hspace{.5in} \downarrow \mbox{\bf 
H}(\varphi)    \\
\mbox{\bf G}(R')    &  \stackrel{f_{R'}}{\longrightarrow}  &    
\mbox{\bf H}(R')
\end{array}
$$
commutes.

Every finite (abstract) group   $\mbox{F}$ gives 
raise to a finite constant $k$-group that will be denoted by 
$\mbox{F}_k$. We have $\mbox{F}_k(R) = \mbox{ F}$ 
whenever $\Spec R$ is connected. Finally, if $\mbox{\bf G}$ is a $k$-group and $R \in k${\it -alg}, we let $\mbox{\bf G}_R$ denote the $R$-group obtained by the base change $k \to R$.

An {\it affine} $k$-group is a $k$-group $\mbox{\bf G}$ 
which is representable: There exists $k[\mbox{\bf G}] \in k${\it 
-alg} such that $\mbox{\bf G}(R) = \Hom(k[\mbox{\bf G}], R)$ for all 
$R$ in $k${\it -alg}. By Yoneda's Lemma $k[\mbox{\bf G}]$ is unique 
up to isomorphism and has a natural commutative Hopf algebra structure. 
Conversely, every commutative Hopf algebra $A$ over $k$ gives rise to an affine 
$k$-group $\mbox{\bf G}^A$ whose $R$-points are $\Hom(A,R)$. 

We 
will make repeated use of the following affine $k$-groups: $\mbox{\bf Hom}\, ( U,V): R \to \Hom_{R\mbox{\it -}mod}(U\otimes R, V\otimes R)$ where $U$ 
and $V$ are finite dimensional $k$-spaces, {\bf End}$\, (V) = {\bf Hom}\, ( V,V)$, 
{\bf G}$_a : R \to (R, +)$, and {\bf G}$_m : R \to R^{\times}$ (the 
units of $R$). In addition we will   make repeated use of  many of the 
classical groups {\bf GL}, {\bf SL}, 
etc; as well as some of the exceptional ones.

Though in the category of affine $k$-groups the concepts of 
kernel of a morphism, normal subgroups, etc. are obvious ones 
(``d\'efinition ensembliste''), those of image of a morphism and 
quotient are not. We illustrate this with an example that will appear 
later on.

Consider the $k$-group {\bf Aut}$\, \gs \gl_n : R \to 
\Aut_{R\mbox{\it -}Lie} (\gs \gl_n(k) \otimes R)$ (this last being the group of  automorphisms 
in the category of Lie algebras over $R$). This is an affine 
$k$-group and there is a natural homomorphism 
$\Ad:\, ${\bf SL}$_n 
\to${\bf Aut}$\,\gs \gl_n$
where $\Ad_R X : A \to XAX^{-1}$ for all 
$A \in \gs \gl_n(k) \otimes R$ and $X \in \mbox{\bf SL}_n(R)$. If $n = 2$, this 
is a {\it surjective} homomorphism, i.e. the image of $\Ad$ equals {\bf 
Aut}$\,\gs \gl_2$. This {\it does not mean} that a typical 
element of $\theta \in ${\bf Aut}$\,\gs \gl_2 (R)$ is of the form $\Ad X$ 
for some $X \in \mbox{\bf SL}_2(R)$ (a false statement even for most 
fields), but rather that there exist a faithfully flat base change 
$R \to S$ and an element $X$ of ${\bf SL}_2(S)$ such that $\theta_S = \Ad_S X$, where 
$\theta_S = \theta \otimes 1 \in \Aut_{S\mbox{\it -}Lie}(\gs \gl_2(k) \otimes R 
\otimes_R S) = ${\bf Aut}$\,\gs \gl_2 (S)$.

The kernel of $\Ad$ is a 
closed subgroup isomorphic to $\boldsymbol{\mu}_n$, where $\boldsymbol{\mu}_n(R):= \{ r \in R 
\, | \, r^n = 1\}$. In particular,  we have the exact sequence $1 \to \boldsymbol{\mu}_n \to 
\mbox{\bf SL}_n \to ${\bf Aut}$\,\gs \gl_n. $  The quotient (in 
the category of affine $k$-groups) $\mbox{\bf SL}_n/ \boldsymbol{\mu}_n$
 is also denoted by $\mbox{\bf PGL}_n$. 
\medskip

We shall make repeated use of the following two standard results

\begin{proposition} \label{morph}
Let $f: \mbox{\bf N} \to \mbox{\bf G}$ and $g: 
\mbox{\bf G} \to \mbox{\bf H}$ be morphisms of affine $k$-groups. The 
following are equivalent.

(i) The sequence $1 \to \mbox{\bf N} 
\stackrel{f}{\rightarrow} \mbox{\bf G} \stackrel{g}{\rightarrow} 
\mbox{\bf H} \to 1$ is exact (i.e. $f$ is injective, $f(\mbox{\bf N}) 
= \ker g$, and $g$ is surjective).

(ii) $f$ is injective and 
$\mbox{\bf G}/f(\mbox{\bf N}) \simeq \mbox{\bf H}$.

(iii) $1 \to 
\mbox{\bf N}(R) \stackrel{f_R}{\rightarrow} \mbox{\bf G}(R) 
\stackrel{g_R}{\rightarrow} \mbox{\bf H}(R)$ is an exact sequence of 
abstract groups and the map $g_k : \mbox{\bf G}(k) \to \mbox{\bf 
H}(k)$ is surjective.

(iv)  $1 \to \mbox{\bf N}(R) 
\stackrel{f_R}{\rightarrow} \mbox{\bf G}(R) 
\stackrel{g_R}{\rightarrow} \mbox{\bf H}(R)$ is an exact sequence of 
abstract groups and for all $h \in \mbox{\bf H}(R)$ there exist a 
faithfully flat base change $R \to S$ and $g \in \mbox{\bf G}(S)$ 
such that $g_S(g) = h_S$ (where $h_S$ is the image of $h$ under the 
map $\mbox{\bf H}(R) \to \mbox{\bf H}(S)$).

\end{proposition}

\noindent
{\bf Proof.} The crucial point is that if $k[\mbox{\bf H}] 
\to k[\mbox{\bf G}]$ is injective then this map is faithfully flat. 
Furthermore, all our $k$-groups are smooth since we are in 
characteristic zero (Cartier's Theorem). 
See \cite{W} and \cite{KMRT}.  \hfill $\square$

\begin{proposition} \label{normal}

Let $\mbox{\bf  G}$ be a $k$-group and 
$\mbox{\bf  H}$ a closed subgroup of  $\mbox{\bf  G}$. The following 
are equivalent:

(i)  $\mbox{\bf  H}$ is normal in  $\mbox{\bf  G}$ 
(i.e.  $\mbox{\bf  H}(R)$ is normal in  $\mbox{\bf  G}(R)$ for all $R 
\in k$-alg).

(ii)  $\mbox{\bf  H}(k)$ is a normal subgroup of  
$\mbox{\bf  G}(k)$.
\end{proposition}
\noindent
{\bf Proof.} Follows 
from \cite{DG}, Ch. II, \S 5.4.1. Indeed, both  $\mbox{\bf  G}$ and  
$\mbox{\bf H}$ are smooth by Cartier's Theorem. \hfill $\square$
\medskip

To illustrate some of these ideas we finish this section by giving the functorial version of a well known result that we shall need later on.

\begin{lemma} \label{schur} (Schur's Lemma) Let $\ga$ be a Lie algebra over $k$, and let $\rho : \ga \to 
\gg \gl (V)$ be an irreducible finite dimensional representation of 
$\ga$. Consider the $k$-group $\mbox{\bf C}_{\rho}$ given by
$$ 
\mbox{\bf C}_{\rho}(R):= \{ \varphi \in \mbox{\bf End}\, V(R)\; | \; 
\varphi \mbox{ commutes with } \rho_R(\ga(R))\}.$$
The canonical map 
{\bf G}$_a \to \mbox{\bf C}_{\rho}$ is an isomorphism. 
\end{lemma}

\noindent
{\bf Proof.} It is clear that $\mbox{\bf C}_{\rho}$ is a 
closed subgroup of {\bf End}$\, V$ and hence affine. If $r \in R$ 
then $r \Id_V \in \mbox{\bf C}_{\rho}(R)$, which leads to a natural 
injection {\bf G}$_a \to \mbox{\bf C}_{\rho}$. By the usual  Schur's 
Lemma  {\bf G}$_a(k)$ surjects onto $\mbox{\bf C}_{\rho}(k)$ and 
therefore our map is an isomorphism of affine $k$-groups. \hfill 
$\square$

\section{The proof of Theorem \ref{th1}} \label{proof1}
\setcounter{equation}{0}

\subsection{Generalities on ${\mathfrak g}{\mathfrak l}(m|n)$ and ${\mathfrak s}{\mathfrak l}(m|n)$  } \label{6.1}

Recall that 
$$
\gg \gl (m|n) := \left\{\left( \begin{array}{cc} A&B\\C&D \end{array} \right) \; | \;
A \in {\cal M}_m, B \in {\cal M}_{m,n}, C \in {\cal M}_{n,m}, D \in {\cal M}_n \right\},
$$
where by ${\cal M}_{p,q} = {\cal M}_{p,q}(k)$ we denote the space of $p\times q$ matrices with entries in $k$ and ${\cal M}_{p}:= {\cal M}_{p,p}$. For $(X,Y) \in  ( \mbox{\bf SL}_m \times \mbox{\bf SL}_n)(R)$ we define $\Ad_R(X,Y) \in \mbox{\bf Aut} \, \gg \gl (m|n)(R)$ by 
$$
\Ad_R(X,Y) : \left( \begin{array}{cc} A&B\\C&D \end{array} \right) \mapsto \left( \begin{array}{cc} XAX^{-1}&XBY^{-1}\\YCX^{-1}&YDY^{-1} \end{array} \right).
$$
This yields a group homomorphism $\Ad : \mbox{\bf SL}_m \times \mbox{\bf SL}_n \to \mbox{\bf Aut} \,\gg \gl (m|n)$. 

The {\it supertranspose} $\tau$ is the automorphism of $\gg \gl (m|n)$ given by 
$$
\tau : \left( \begin{array}{cc} A&B\\C&D \end{array} \right) \to \left( \begin{array}{cc} -A^{t}&C^{t}\\-B^{t}&-D^{t} \end{array} \right).
$$
The automorphism $\tau$ is of order $4$ and leads to a constant subgroup $<\tau>_k$ of $\mbox{\bf Aut} \,\gg \gl (m|n)$. 

If $m=n$ we have an automorphism $\pi$ of $\gg \gl (n|n)$ given by $$\pi: \left( \begin{array}{cc} A&B\\C&D \end{array} \right) \to \left( \begin{array}{cc} D&C\\B&A \end{array} \right)$$
which leads to a constant group $<\pi>_k$ of $\mbox{\bf Aut} \,\gg \gl (n|n)$. 

We also have an injective group homomorphism $j: \mbox{\bf G}_m \to \mbox{\bf Aut}\, \gg \gl (m|n)$ given by $j_R(\lambda): \left( \begin{array}{cc} A&B\\C&D \end{array} \right) \to \left( \begin{array}{cc} A&\lambda B\\\lambda^{-1}C&D \end{array} \right)$ for all $\lambda \in  \mbox{\bf G}_m(R) = R^{\times}$.

One easily checks that:
\begin{equation} \label{iden}
\begin{array}{c} 
\tau_R \Ad_R(X,Y) \tau_R^{-1} = \Ad_R((X^{t})^{-1},(Y^{t})^{-1}), \pi_R \Ad_R(X,Y) \pi_R^{-1} = \Ad_R(Y, X), \\
\pi_R \tau_R \pi_R^{-1} = \tau_R^{3}, \tau_R^{2} = j_R(-1), \tau_R j_R(\lambda) \tau_R^{-1} = j_R(\lambda^{-1}).
\end{array}
\end{equation}

\subsection{Automorphisms of ${\mathfrak s}{\mathfrak l} (m|n), m>n\geq 1$}

Recall that $\gs \gl (m|n)$ consists of those elements $\left( \begin{array}{cc} A&B\\C&D \end{array} \right)$ of $\gg \gl (m|n)$ for which $\tr (A) - \tr(D) = 0$. As above we have a group morphism $\Ad \times j : \mbox{\bf SL}_m \times \mbox{\bf SL}_n \times \mbox{\bf G}_m \to \mbox{\bf Aut}\, \gs \gl (m|n)$. The image of $\Ad \times j$ is denoted by $\mbox{\bf Aut}^0 \gs \gl (m|n)$. Clearly $(X, Y, \lambda) \in \ker (\Ad \times j)_R$ iff there exist $\alpha, \beta \in R^{\times}$ such that $X = \alpha I_m, Y = \beta I_n, \alpha^m = \beta^n = 1$, and $\lambda = \alpha^{-1} \beta$. Thus $\ker (\Ad \times j) = \boldsymbol{\mu}_m \times \boldsymbol{\mu}_n$.

Serganova's reasoning shows that $\mbox{\bf Aut}\, \gs \gl (m|n)(k)$ is generated by $\mbox{\bf Aut}^0 \gs \gl (m|n)(k)$ and $\tau$. From Proposition \ref{morph} and \ref{normal} it follows that $\mbox{\bf Aut}^0 \gs \gl (m|n)$ is normal in  $\mbox{\bf Aut}\, \gs \gl (m|n)$ and that the canonical morphism $<\tau>_k \to \mbox{\bf Aut}\, \gs \gl (m|n)/ \mbox{\bf Aut}^0 \gs \gl (m|n)$ is surjective. With the aid of (\ref{iden}) and Proposition \ref{morph} the kernel of this last homomorphism is easily seen to be $<\tau^2>_k$. 

The above definition yields a non-split exact sequence $$1 \to (\mbox{\bf SL}_m \times \mbox{\bf SL}_n \times \mbox{\bf G}_m)/ (\boldsymbol{\mu}_m \times \boldsymbol{\mu}_n) \to \mbox{\bf Aut} \, \gs \gl (m|n) \to \Z_{2,k} \to 1,$$ where $<\tau>_k$ surjects onto $\Z_{2,k}$.

\subsection{Automorphisms of ${\mathfrak p}{ \mathfrak s}{\mathfrak l} (n|n), n > 2$}

Recall that  $\gp \gs \gl (n|n) := \gs \gl (n|n) / \gz$ where $\gz:=k \left( \begin{array}{cc} I_n&0\\0&  I_n \end{array} \right)$. One reasons as in the last case to obtain a surjective group homomorphism $<\tau>_k \rtimes <\pi>_k \to \mbox{\bf Aut} \gp \gs \gl (m|n)/ \mbox{\bf Aut}^0 \gp \gs \gl (m|n)$ (see Proposition \ref{morph}). An easy reasoning again shows that the kernel of this map is $<\tau^2>_k$. Since the corresponding quotient is the Klein $4$-group we have $$1 \to (\mbox{\bf SL}_n \times \mbox{\bf SL}_n \times \mbox{\bf G}_m)/ (\boldsymbol{\mu}_n \times \boldsymbol{\mu}_n) \to \mbox{\bf Aut} \, \gp \gs \gl (n|n) \to \Z_{2,k} \times \Z_{2,k} \to 1.$$
The sequence is not split, but the second copy of $\Z_{2,k}$ (which corresponds to $<\pi>_k$) does admit a section.

\subsection{Automorphisms of ${\mathfrak p}{\mathfrak s}{\mathfrak l}(2|2)$}

In this case $\gg_{\bar{1}} = \gg_{-1} \oplus \gg_1$, $\gg_{-1} \simeq  V_2^* \otimes V_2$, and $\gg_1 \simeq V_2 \otimes V_2^*$, where $V_2$ is the standard $\gg \gl_2$-module. In fact, $\gg_{-1}$ and $\gg_1$ are isomorphic $\gg_0$-modules. We explicitly describe how this is done to better understand Serganova's argument. Let $\{ e_1, e_2\}$ be the standard basis of $V_2$, and $\{ e_1^*, e_2^*\}$ the corresponding dual basis of $V_2^*$. Then $\varphi: V_2 \stackrel{\simeq}{\rightarrow} V_2^*$ as $\gs \gl_2$-modules via $\varphi(e_1) = - e_2^*, \varphi(e_2) =  e_1^*$. This yields the isomorphism $\phi:= \varphi \times \varphi^{-1}: V_2 \otimes V_2^* \to V_2^* \otimes V_2$. The matrical description of $\phi$ is as follows. Identify $V_2 \otimes V_2^*  \stackrel{\eta}\simeq  {\cal M}_2(k)$ via $e_i \otimes e_j^* \mapsto E_{ij}$, and  $V_2^* \otimes V_2  \stackrel{\eta^*}\simeq {\cal M}_2(k)$ via $e_i^* \otimes e_j \mapsto E_{ij}$. We then have the following commutative diagram:

$$
\begin{array}{ccc}
V_2 \otimes V_2^*  &  \stackrel{\phi}{\longrightarrow} &  V_2^* \otimes V_2 \\
\eta \downarrow   &       &  \downarrow \eta^*  \\
{\cal M}_2(k)   &  \stackrel{\psi}{\longrightarrow}    &    {\cal M}_2(k)
\end{array}
$$
where $\psi(E) = -JE^tJ^{-1}$ with $J = \left( \begin{array}{cc} 0&1\\-1&0 \end{array} \right)$. Furthermore, $\psi([B,C]) = [\psi(B), \psi(C)]$ for all $B$ and $C$ in ${\cal M}_2(k)$. Using the canonical $R$-module isomorphism $(V_2\otimes R)^* \simeq V_2^* \otimes R$, one easily sees that the above construction is ``functorial on $R$''. More precisely, we have the commutative diagram

\begin{equation} \label{psl2}
\begin{array}{cccccccc}
\phi_R:&  V_2(R) \otimes V_2(R)^* & \simeq  & (V_2 \otimes V_2^*)(R)&  \stackrel{\phi \otimes 1}{\longrightarrow} &  (V_2^* \otimes V_2)(R) & \simeq  & V_2(R)^* \otimes V_2(R)\\
 & \eta_R \downarrow \hspace{.3in}  &  & \eta \otimes 1 \downarrow \hspace{.3in} &     & \hspace{.3in} \downarrow \eta^* \otimes 1 &  & \hspace{.2in} \downarrow \eta_R^*  \\
\psi_R: & {\cal M}_2(R) & \simeq  & {\cal M}_2(k) \otimes R &   \stackrel{\psi \otimes 1}{\longrightarrow} & {\cal M}_2(k) \otimes R &\simeq    &    {\cal M}_2(R)
\end{array}
\end{equation}
where the compositions $\phi_R$  and $\psi_R$  of all horizontal isomorphisms are  $(\gs \gl_2 \oplus \gs \gl_2) (R)$-module isomorphisms.

For $\left( \begin{array}{cc} \alpha&\beta\\ \gamma&\delta \end{array} \right) \in \mbox{\bf SL}_2(R)$ consider the $R$-module endomorphism 
$\rho_R \left( \begin{array}{cc} \alpha&\beta\\ \gamma&\delta \end{array} \right)$ of $\gp \gs \gl (2|2)$ given by 
$$
\rho_R \left( \begin{array}{cc} \alpha&\beta\\ \gamma&\delta \end{array} \right) : \left( \begin{array}{cc} A & B\\ C & D \end{array} \right) \mapsto \left( \begin{array}{cc} A &\alpha B + \beta \psi_R(C)\\ \gamma \psi_R(B) + \delta C& D \end{array} \right).
$$
A direct calculation using (\ref{psl2}) shows that this is in fact an element of $\mbox{\bf Aut} \, \gp \gs \gl (2|2)$, and that this construction yields a $k$-group homomorphism $\rho : \mbox{\bf SL}_2 \to \mbox{\bf Aut} \, \gp \gs \gl (2|2) $ with trivial kernel. A direct computation shows that $\rho(\mbox{\bf SL}_2 )$ commutes with the image of $\Ad$. This leads to a group homomorphism $$\Ad \times \rho : \mbox{\bf SL}_2 \times \mbox{\bf SL}_2 \times \mbox{\bf SL}_2 \to \mbox{\bf Aut} \, \gp \gs \gl (2|2) $$ with kernel $\boldsymbol{\mu}_2 \times \boldsymbol{\mu}_2$ corresponding to $(-I_2, I_2, -I_2)$ and $(I_2, -I_2, -I_2)$. Furthermore, $\pi_R (\Ad \times \rho)_R (X,Y,Z) \pi_R^{-1} = (\Ad \times \rho)_R(Y,X,(Z^t)^{-1})$ and it easily follows from Proposition \ref{normal} that the image of $\Ad \times \rho$ is normal in $\mbox{\bf Aut} \, \gp \gs \gl (2|2) $. Serganova's argument shows that $\mbox{\bf Aut} \, \gp \gs \gl (2|2)(k)$ is generated by the image of $(\Ad \times \rho )_k$ and $\pi$. By Proposition \ref{morph}, we obtain the split exact sequence
$$
1 \to (\mbox{\bf SL}_2 \times \mbox{\bf SL}_2 \times\mbox{\bf SL}_2 )/(\boldsymbol{\mu}_2 \times \boldsymbol{\mu}_2) \to \mbox{\bf Aut} \, \gp \gs \gl (2|2) \to \Z_{2, k} \to 1.
$$

\subsection{Automorphisms of $\mathfrak g = \gs \pp(n)$}

The Lie subsuperalgebra $\pp(n)$ of $\gg \gl (n|n)$ in matrix form is defined by $$\pp(n):=\left\{ \left( \begin{array}{cc} A& B \\ C & -A^t  \end{array} \right)\; | \; A \in \gg \gl_n, B = B^t, C=-C^t \right\}.$$ By $\gs \pp (n)$ we denote the quotient $\pp(n)/\gz$ where $\gz:= \left\{ \left(\begin{array}{cc} \lambda I_n& 0 \\ 0 & -\lambda I_n  \end{array} \right)\; | \; \lambda \in k \right\}$. We have that $\gs \pp(n)_{\bar{0}} \simeq \gs \gl_n$, $\gs \pp (n)_{\bar{1}} = \gs \pp (n)_{1} \oplus \gs \pp (n)_{ -1}$, $ \gs \pp (n)_{1} \simeq S^2 V_n$, and  $\gs \pp (n)_{-1} \simeq \Lambda^2 V_n$.

There is a natural homomorphism $\Ad: \mbox{\bf SL}_n  \to \mbox{\bf Aut} \, \gg$ given by $X \mapsto \Ad (X, X^*)$, where $X^* := (X^t)^{-1}$. Explicitly $$\Ad_R(X) : \left( \begin{array}{cc} A& B \\ C & D  \end{array} \right) \mapsto \left( \begin{array}{cc} XAX^{-1}& XBX^t \\ X^*CX^{-1} & X^*DX^t  \end{array} \right)$$ for all $X \in \mbox{\bf SL}_n(R)$ and $\left( \begin{array}{cc} A& B \\ C & D  \end{array} \right) \in \gs \pp(n)(R)$. Combined with the map $j$ as above this yields a group homomorphism 
$$
\Ad \times j : \mbox{\bf SL}_n \times  \mbox{\bf G}_m \to \mbox{\bf Aut} \,  \gs \pp(n)
$$
that can be seen to be surjective on $k$-points by reasoning as Serganova's. Clearly $(X, \lambda)$ belongs to $\ker (\Ad \times j)_R$ iff $X = \alpha I_n$ with $\alpha^n = 1$ and $\alpha^2 \lambda = 1$. Thus:
$$
\mbox{\bf Aut} \,  \gs \pp(n) \simeq (\mbox{\bf SL}_n \times  \mbox{\bf G}_m )/ \boldsymbol{\mu}_n.
$$

\subsection{Automorphisms of ${\mathfrak p}{\mathfrak s}{\mathfrak q}(n)$}

The Lie subsuperalgebra $\gq(n)$ of $\gg \gl (n|n)$ is defined in matrix form by 
$$
\gq(n):= \left\{ \left( \begin{array}{cc} A& B \\ B & A  \end{array} \right)\; | \; A,B \in \gg \gl_n    \right\}.
$$
We denote by $\gs \gq (n)$ the subsuperalgebra of odd-traceless matrices $\left\{ \left( \begin{array}{cc} A& B \\ B & A  \end{array} \right)\; | \; \tr (B)=0 \right\}$, and by $\gp \gs \gq (n)$ the quotient $\gs \gq (n) / \gz$, where $\gz := \left\{\left( \begin{array}{cc} \lambda I_n & 0 \\ 0 & \lambda I_n  \end{array} \right)\; | \; \lambda \in k \right\}$. We have also that $\gp \gs \gq (n)_{\bar{0}} \simeq \gs \gl_n$ and $\gp \gs \gq (n)_{\bar{1}} \simeq \gs \gl_n$ is isomorphic to the adjoint representation of $\gs \gl_n$.  Note that  $\gq(n)$ is not invariant with respect to the supertransposition, however it is so with respect to the $q$-{\it supertransposition} $\sigma_q : \left( \begin{array}{cc} A& B \\ B & A  \end{array} \right)\mapsto  \left( \begin{array}{cc} A^t& \zeta B^t \\ \zeta B^t & A^t  \end{array} \right)$, where $\zeta = \zeta_4 \in k$ is our fixed primitive $4$th root of unity.	
There is a natural homomorphism $\mbox{\bf SL}_n \to \mbox{\bf Aut} \,  \gp \gs \gq (n)$ arising from $X \mapsto \Ad_R(X,X)$. The kernel of this map is $\boldsymbol{\mu}_n$ showing that $\mbox{\bf Aut} \,  \gp \gs \gq (n)$ contains a copy of $\mbox{\bf PGL}_n$. We also have a constant subgroup $<\sigma_q>_k \simeq \Z_{4,k}$ of  $\mbox{\bf Aut} \,  \gp \gs \gq (n)$. The considerations of \S \ref{Gen} and Serganova's original argument show that we have a split exact sequence
$$
1 \to \mbox{\bf PGL}_n \to \mbox{\bf Aut} \,  \gp \gs \gq (n) \to \Z_{4,k} \to 1.
$$

\subsection{Automorphisms of ${\mathfrak o}{\mathfrak s}{\mathfrak p}(m|2n)$} \label{osp}

The orthosymplectic Lie superalgebra is defined by 
$$ 
\go \gs \gp(m|2n) := \{ Z \in \gg \gl (m|2n) \; | \; Z^{st}B_{m,n} + B_{m,n}Z =0\},
$$
where $B_{m,n}:= \left( \begin{array}{cc} I_m& 0 \\ 0 & J_{n}  \end{array} \right)$ and $J_n:= \left( \begin{array}{cc} 0 & I_n \\ -I_n & 0  \end{array} \right)$.  Equivalently, we have 
$$\go \gs \gp(m|2n)  = \left\{ \left( \begin{array}{cc} A& B \\ J_{n}B^t & D  \end{array} \right)\; | \; A \in \gs\go_m, B \in {\cal M}_{m,2n}, D \in \gs \gp_{2n}\right\}.
$$
We have that  $\go \gs \gp(m|2n)_{\bar{0}} = \gs \go_m \oplus \gs \gp_{2n}$ and  $\go \gs \gp(m|2n)_{\bar{1}}$ is the $ (\gs \go_m \oplus \gs \gp_{2n})$-module $V_m \otimes V_{2n}$. Note that $\go \gs \gp(2|2n)_{\bar{1}} \simeq (k\otimes V_{2n}) \oplus (k\otimes V_{2n})$.

We have a natural homomorphism 
$$
\Ad : \mbox{\bf O}_m \times \mbox{\bf Sp}_{2n} \to \mbox{\bf Aut} \, \go \gs \gp (m|2n)
$$
which is surjective on $k$-points as shown by Serganova. The kernel of $\Ad$ is isomorphic to $\boldsymbol{\mu}_2$ coming from $\Ad(-I_m, -I_{2n})$. Thus 
$$
\mbox{\bf Aut} \, \go \gs \gp (m|2n) \simeq (\mbox{\bf O}_m \times \mbox{\bf Sp}_{2n}) / \boldsymbol{\mu}_2.
$$
Consider the canonical map $\tau: \mbox{\bf SO}_m \times \mbox{\bf Sp}_{2n} \to (\mbox{\bf O}_m \times \mbox{\bf Sp}_{2n}) / \boldsymbol{\mu}_2$. If $m$ is odd then $\tau$ is injective and  clearly surjective on $k$-points because $(X, Y) \cong ((-I_m)X, (-I_{2n})Y) \mod \boldsymbol{\mu}_2(k)$. Thus 
$$
\mbox{\bf Aut} \, \go \gs \gp (m|2n) \simeq \mbox{\bf SO}_m \times \mbox{\bf Sp}_{2n} \mbox{ if $m$ is odd}.
$$
Assume now that $m$ is even. Let $r_m \in \mbox{\bf O}_m(k)$ be an element satisfying $\det r_m = -1$ and $r_m^2 = I_m$. Let $r:=\Ad(r_m, I_{2n}) \in \mbox{\bf Aut} \, \go \gs \gp (m|2n)$. Then $r$ leads to a constant subgroup $<r>_k \simeq \Z_{2,k}$. It is immediate that $r_{R} \Ad_R(X, Y) r_{R}^{-1} = \Ad_R(r_{m, R} Xr_{m, R}^{-1}, Y)$ where $r_{m, R}$ is the image of $r_m$ under the canonical map $\mbox{\bf O}_m(k) \to \mbox{\bf O}_m(R)$. This yields a group homomorphism 
$$
(\mbox{\bf SO}_m \times \mbox{\bf Sp}_{2n})\rtimes <r>_k \to \mbox{\bf O}_m \times \mbox{\bf Sp}_{2n} \to  (\mbox{\bf O}_m \times \mbox{\bf Sp}_{2n})/\boldsymbol{\mu}_2 \simeq  \mbox{\bf Aut} \, \go \gs \gp (m|2n) 
$$
which is surjective on $k$-points by Serganova's original reasoning. Since the kernel is precisely the group $\boldsymbol{\mu}_2 \times I$ described above we have the split exact sequence
$$
1 \to (\mbox{\bf SO}_m \times \mbox{\bf Sp}_{2n})/\boldsymbol{\mu}_2 \to \mbox{\bf Aut} \, \go \gs \gp (m|2n) \to \Z_{2,k} \to 1
$$
whenever $m$ is even.

\subsection{The exceptional Lie superalgebras} \label{except}
 In this subsection we focus on the three exceptional Lie superalgebras $\gg = G(3)$, $\gg = F(4)$, and $\gg = D(\alpha)$. First we recall some basic facts about the multiplication structure of each $\gg$. For more details we refer the reader to \cite{K2} and \cite{S}. To describe the Lie bracket on $\gg$ we introduce some notations. Denote by $p$ the bilinear form on $V_2$ with matrix $J$, i.e. $p(u,v) := u^{t}Jv$ for $u, v \in V_2$. For $c \in k^{\times}$ we define a map $\Psi^c:V_2 \otimes V_2 \to \gs \gl_2$ by the formula $\Psi^c(u,v)w:= c(p(v,w)u - p(w,u)v)$ and  write simply $\Psi$ for $\Psi^1$. One checks immediately that if $A \in \mbox{\bf GL}_2(k)$ then  $A\Psi^c(u,v)A^{-1} = \Psi^c(Au,Av)$ for all $u,v$ iff $p(Au,Av) = p(u,v)$ for all $u,v$ iff $\det A = 1$.

\medskip
\noindent
{\bf \ref{except}.1 The case $\gg = G(3)$.} We have $\gg_{\bar{0}} \simeq G_2 \oplus \gs \gl_2$ and $\gg_{\bar{1}}$ is isomorphic to the simple $\gg_{\bar{0}}$-module $stan_{G_2} \otimes V_2$, where $G_2$ is the exceptional Lie algebra $\Der({\mathbb O})$ of the derivations of the octonions $\mathbb O$, and  $stan_{G_2}$ is the standard $8$-dimensional $G_2$-module $\mathbb O$. Let $\pi_1: G_2 \oplus \gs \gl_2 \to G_2$ and $\pi_2: G_2 \oplus \gs \gl_2 \to \gs \gl_2$ be the natural projections.   The projection $\pi_2$ of the restriction $[\cdot,\cdot]_{|\gg_{\bar{1}} \times \gg_{\bar{1}}}: \gg_{\bar{1}} \times \gg_{\bar{1}} \to \gg_{\bar{0}}$ of the Lie bracket on $\gg_{\bar{1}}$ is described by the formula: $\pi_2 [x_1\otimes y_1, x_2 \otimes y_2] = (x_1,x_2)_{\mathbb O}\Psi(y_1,y_2)$, where $x_i \otimes y_i \in stan_{G_2} \otimes V_2 \simeq \gg_{\bar{1}} $ and  $(x_1,x_2)_{\mathbb O}:=x_1x_2^* + x_2 x_1^*$ is the non-degenerate $G_2$-invariant form on $\mathbb O$ (for more details see for example \cite{F}). The formula of $\pi_1$, on the other hand, comes from the adjoint representation.

It is clear then that we have an injective $k$-group homomorphism $\Ad : \mbox{\bf G}_2 \times \mbox{\bf SL}_2 \to \mbox{\bf Aut}\, G(3)$ given by:
$$
\Ad_R(X,Y) : ((G,A), u\otimes v) \mapsto ((XGX^{-1}, YAY^{-1}), Xu \otimes Yv).
$$
Since by Serganova's original work this map is surjective on $k$-points we obtain:
$$
\mbox{\bf Aut}\, G(3) \simeq  \mbox{\bf G}_2 \times \mbox{\bf SL}_2.
$$

\medskip
\noindent
{\bf \ref{except}.2 The case $\gg = F(4)$.} We have $\gg_{\bar{0}} \simeq \gs \go_7 \oplus \gs \gl_2$ and $\gg_{\bar{1}}$ is isomorphic to the simple $\gg_{\bar{0}}$-module $spin_{7} \otimes V_2$, where $spin_7$ is the $\Spin$-representation of $\gs \go_7$ whose space is the Clifford algebra $\Cliff_6$. If $\pi_1$ and $\pi_2$ are the canonical projections of $\gg_{\bar{0}}$ to its first and second summands, respectively, then we have again that $\pi_2 [x_1\otimes y_1, x_2 \otimes y_2] = (x_1,x_2)\Psi(y_1,y_2)$ where $x_i \otimes y_i \in spin_7 \otimes V_2 \simeq \gg_{\bar{1}} $ and $(\cdot,\cdot)$ is an $\gs \go_7$-invariant form on $spin_7$. Similarly to the case $\gg = G(3)$ we conclude that 
$$
\mbox{\bf Aut}\,F(4) \simeq (\mbox{\bf Spin}_7 \times \mbox{\bf SL}_2)/{\boldsymbol \mu_2}.
$$
Here $\mbox{\bf Spin}_7$ is the simply connected $k$-group of type $B_3$ and ${\boldsymbol \mu_2}$ is the diagonal subgroup of the center ${\boldsymbol \mu_2} \times {\boldsymbol \mu_2}$ of $\mbox{\bf Spin}_7 \times \mbox{\bf SL}_2$.

\medskip
\noindent
{\bf \ref{except}.3 The case $\gg = D(\alpha)$} (often denoted by $D(2,1,\alpha)$). Now we have $\gg_{\bar{0}} \simeq \gs \gl_2 \oplus \gs \gl_2 \oplus \gs \gl_2$ and $\gg_{\bar{1}} \simeq V_2\otimes V_2 \otimes V_2$. Denote by $(\gs \gl_2)_i$ and $(V_2)_i$ the $i$-th components of $\gg_{\bar{0}}$ and $\gg_{\bar{1}}$ respectively and by $\pi_i : \gg_{\bar{0}} \to (\gs \gl_2)_i$ the natural projections.  Then for $\alpha_1:= \alpha, \alpha_2 := 1,\alpha_3 := -1 - \alpha$ and $u_i \in (V_2)_i$ we have 
$$
[u_1\otimes u_2 \otimes u_3 , v_1 \otimes v_2 \otimes v_3] = 
\sum_{\sigma \in {\mathfrak S}_3}p_{\sigma(1)}(u_{\sigma(1)}, v_{\sigma(1)})p_{\sigma(2)}(u_{\sigma(2)}, v_{\sigma(2)})\Psi_{\sigma(3)}^{\alpha_{\sigma(3)}}(u_{\sigma(3)}, v_{\sigma(3)}),
$$
 where $p_j$ and $\Psi_j^c$ are the maps $p$ and $\Psi^c$ on $(\gs \gl_2)_j$. Note that $D(\alpha) \simeq D(\alpha')$ if and only if $\alpha' = -(\alpha+1)^{\pm 1}$,  $\alpha' = \alpha^{-1}$, or  $\alpha' = -(\alpha/(\alpha +1))^{\pm 1}$. In future we consider $\alpha \neq 0, -1$ since $D(0)$ and $D(-1)$ are not simple.

The natural homomorphism $\Ad : \mbox{\bf SL}_2 \times\mbox{\bf SL}_2  \times \mbox{\bf SL}_2 \to  \mbox{\bf Aut} \,  D(\alpha)$ has kernel $\boldsymbol{\mu}_2 \times \boldsymbol{\mu}_2$ coming from $(-I_2, I_2, - I_2)$ and $(I_2, -I_2, -I_2)$. As usual we let $\mbox{\bf Aut}^0   D(\alpha)$ denote the image of $\Ad$. Because $\mbox{\bf Aut}\, (\gs \gl_2 \oplus \gs \gl_2 \oplus \gs \gl_2 ) \simeq (\mbox{\bf PGL}_2)^3 \rtimes {\mathfrak S}_{3, k}$  the outer group can be found inside ${\mathfrak S}_{3, k}$. The answer depends on the value of $\alpha$. Fix $\sigma \in {\mathfrak S}_3$ and $\lambda \in k^{\times}$. Define $\theta (\sigma, \lambda) \in \mbox{\bf GL}(D(\alpha))$ by 
$$
\theta (\sigma, \lambda) : ((A_1, A_2, A_3), u_1 \otimes u_2 \otimes u_3) \to ((A_{\sigma(1)}, A_{\sigma(2)}, A_{\sigma(3)}), \lambda  u_{\sigma(1)} \otimes u_{\sigma(2)} \otimes u_{\sigma(3)}).
$$

Following Serganova we verify that:

(i) If $\alpha =1$ then $\theta ((1,2), 1) \in \mbox{\bf Aut} \,  D(\alpha)$,

(ii) If  $\alpha^3 =1, \alpha \neq 1$ then $\theta ((1,2,3), \lambda) \in \mbox{\bf Aut} \,  D(\alpha)$, where $\lambda^2 = \frac{1}{\alpha}$.

A straightforward calculation shows that 
$$\theta(\sigma, \lambda)_R \Ad(X_1, X_2, X_3) \theta (\sigma, \lambda)_R^{-1} = \Ad (X_{\sigma(1)}, X_{\sigma(2)},X_{\sigma(3)}).$$
We thus have the following injective $R$-group homomorphisms:
$$
\begin{array}{rccl}

\mbox{\bf Aut}^0 D(\alpha) & \to & \mbox{\bf Aut} \, D(\alpha)  &  \mbox{ for }\alpha \notin \{1, -1/2, -2\}\mbox{ and }\alpha^3 \neq 1,\\

\mbox{\bf Aut}^0 D(1) \rtimes <\theta((1,2), 1)>_k & \to & \mbox{\bf Aut}\, D(1)  & \\

\mbox{\bf Aut}^0 D(\alpha) \rtimes <\theta((1,2,3), \lambda)>_k & \to & \mbox{\bf Aut} \, D(\alpha)  & \mbox{ for }\alpha \neq 1 \mbox{ and }\alpha^3 =  1.

\end{array}
$$
Serganova's original argument shows that these three homomorphisms are surjective on $k$-points. Since $D(\alpha) \simeq D(\beta)$ iff $\alpha, \beta \in \{ 1, -1/2, -2\}$, or $\alpha^3 = \beta^3 = 1$ and $\alpha \neq 1, \beta \neq 1$, we obtain:

(a) If $\alpha \notin \{ 1, -1/2, -2\}$ and $\alpha^3 \neq 1$ then
$$
\mbox{\bf Aut}\, D(\alpha) \simeq (\mbox{\bf SL}_2 \times \mbox{\bf SL}_2 \times \mbox{\bf SL}_2) / (\boldsymbol{\mu}_2 \times \boldsymbol{\mu}_2).
$$

(b) If $\alpha \in \{ 1, -1/2, -2\}$ there is a split exact sequence
$$
1 \to  (\mbox{\bf SL}_2 \times \mbox{\bf SL}_2 \times \mbox{\bf SL}_2)  / (\boldsymbol{\mu}_2 \times \boldsymbol{\mu}_2) \to \mbox{\bf Aut}\, D(\alpha) \to \Z_{2,k} \to 1.
$$

(c) If $\alpha^3 = 1$ and $\alpha \neq 1$ there is a split exact sequence 
$$
1 \to  (\mbox{\bf SL}_2 \times \mbox{\bf SL}_2 \times \mbox{\bf SL}_2)  / (\boldsymbol{\mu}_2 \times \boldsymbol{\mu}_2) \to \mbox{\bf Aut}\, D(\alpha) \to \Z_{3	,k} \to 1.
$$

\subsection{The Cartan type Lie superalgebras} \label{cartan}

In this subsection we consider the simple Cartan type Lie superalgebras $W(n)$, $S(n)$, $S'(n)$ (for $n=2k$), and $H(n)$. A brief summary for these Lie superalgebras is as follows (see \cite{K2} for more details). First recall that every $\Z$-graded associative algebra $A = \oplus_{i \in \Z} A_i$ has a natural $\Z_2$-grading $A = A_{\bar{0}} \oplus A_{\bar{1}}$ where $ A_{\bar{0}}:=\oplus_{i \in \Z}A_{2i}$ and   $ A_{\bar{1}}:=\oplus_{i \in \Z} A_{2i+1}$. $A$ becomes a Lie superalgebra defining a Lie superbracket on $A$ by the equality $[a,b] := ab - (-1)^{\deg a \deg b}ba$, where $\deg x = \bar{i}$ whenever $x \in A_{\bar{i}}$.  (We will also write $\deg x = i$ for  $x \in A_i$). By $W(n)$ we denote the (super)derivations of the Grassmann algebra $\Lambda(n) :=\Lambda(\xi_1,...,\xi_n)$. Every  element $D$ of $W(n)$ has the form $D = \sum_{i=1}^n P_i(\xi_1,...,\xi_n) \frac{\partial}{\partial \xi_i}$ where by definition $\frac{\partial}{\partial \xi_i} (\xi_j) = \delta_{ij}$. Both $\Lambda(n)$ and $W(n)$ have natural gradings $\Lambda(n) = \oplus_{i=0}^n \Lambda(n)_i$ and  $W(n) = \oplus_{j=-1}^{n-1} W(n)_j$, where $\Lambda(n)_i:= \{ P(\xi_1,...,\xi_n) \; | \; \deg P = i\}$ and  $W(n)_{j}:=\{\sum_{i=1}^n P_i \frac{\partial}{\partial \xi_i} \; | \; \deg P_i = j+1\}$.  Note that $[\frac{\partial}{\partial \xi_i}, \frac{\partial}{\partial \xi_j}] = \frac{\partial}{\partial \xi_i} \frac{\partial}{\partial \xi_j} + \frac{\partial}{\partial \xi_j} \frac{\partial}{\partial \xi_i} =0 $ (the first equality holds because $\deg \frac{\partial}{\partial \xi_i} = \bar{1}$).

There are two associative superalgebras of differential forms defined over $\Lambda(n)$, namely $\Omega(n)$ and $\Theta(n)$. The superalgebra $\Omega(n)$ has generators $d \xi_1,...,d \xi_n$ and defining relations $d \xi_i \circ d \xi_j = d \xi_j  \circ d \xi_i$, $\deg d\xi_i = \bar{0}$, while the superalgebra $\Theta(n)$ has generators  $\theta \xi_1,...,\theta \xi_n$ and relations $\theta  \xi_i \wedge \theta \xi_j = - \theta \xi_j  \wedge \theta \xi_i, \deg \theta \xi_i = \bar{1}$. Note that the differentials $d$ and $\theta$ are  derivations of degrees $\bar{1}$ and  $\bar{0}$ respectively. Let $\Delta_n:= \theta \xi_1 \wedge ... \wedge \theta \xi_n$ be the standard volume form in $\Omega(n)$, and let  $\omega_n := \sum_{i=1}^n d \xi_i \circ d \xi_i$ be the standard Hamiltonian form in $\Theta(n)$. Put $\Delta_n':=(1 + \xi_1 ...\xi_n)\Delta_n$. 

Every derivation $D$ of $W(n)$ and every automorphism $\Phi$ of $\Lambda(n)$ extend uniquely to a derivation $\widetilde{D}^{d}$ and an automorphism $\widetilde{\Phi}^{d}$ (respectively, $\widetilde{D}^{\theta}$ and $\widetilde{\Phi}^{\theta}$) of $\Omega(n)$ (resp., $\Theta(n)$) so that $[\widetilde{D}^{d}, d]= 0$ and $[\widetilde{\Phi}^{d}, d]= 0$  (resp., $\widetilde{D}^{\theta} \theta f - \theta  \widetilde{D}^{\theta}f =0$ and  $\widetilde{\Phi}^{\theta} \theta f - \theta  \widetilde{\Phi}^{\theta}f =0$  for every $f \in \Lambda(n)$). We denote by $S(n)$ the Lie superalgebra $\{ D \in W(n)\; | \; \widetilde{D}^\theta(\Delta_n) = 0\} $ and by $S'(n)$ (not to be confused with the derived superalgebra $S(n)'$ of $S(n)$) the Lie superalgebra  $\{ D \in W(n)\; | \; \widetilde{D}^\theta(\Delta_n') = 0\} $. Since $S'(2k+1) \simeq S(2k+1)$ we are interested in $S'(n)$ only for even numbers $n$. The Hamiltonian Lie superalgebras are defined by  $\widetilde{H}(n):= \{ D \in W(n)\; | \; \widetilde{D}^d\omega_n = 0\}$ and $H(n):=[\widetilde{H}(n), \widetilde{H}(n)]$. The following explicit description of the Cartan type series $S$, $S'$, $H$, and $\widetilde{H}$ is very helpful for our considerations: 
$$
S(n) = \Span_k \left\{ \frac{\partial f}{\partial \xi_i} \frac{\partial}{\partial \xi_j} +  \frac{\partial f}{\partial \xi_j} \frac{\partial}{\partial \xi_j}\; | \; f \in \Lambda(n), i,j = 1,...,n\right\},
$$
$$
S'(n) = \Span_k \left\{ (1 - \xi_1 ... \xi_n)\left(\frac{\partial f}{\partial \xi_i} \frac{\partial}{\partial \xi_j} +  \frac{\partial f}{\partial \xi_j} \frac{\partial}{\partial \xi_j}\right)\; | \; f \in \Lambda(n), i,j = 1,...,n\right\},
$$
$$
\widetilde{H}(n) = \Span_k \left\{D_f:=\sum_i \frac{\partial f}{\partial \xi_i} \frac{\partial}{\partial \xi_i} \; | \; f \in \Lambda(n), f(0)=0, i,j = 1,...,n\right\},
$$
$$
\widetilde{H}(n) = H(n) \oplus kD_{\xi_1...\xi_n}.
$$ 
We have also that $[D_f,D_g] = D_{\{f,g\}}$ where $\{f,g\}:=(-1)^{\deg f}\sum_{i=1}^n \frac{\partial f}{\partial \xi_i} \frac{\partial g}{\partial \xi_j}$.

For any Lie subsuperalgebra $\gk$ of $W(n)$ define $\gk_i:= \{ D \in \gk\; | \; \deg D = i\}$ and $\gk^i:= \oplus_{j \geq i} \gk_j$.  If $\gk$ is any of the following Lie superalgebras $W(n), S(n), S'(n), \widetilde{H}(n), H(n)$ then $\gk = \gk^{-1} \supset \gk^{0} \supset...\supset \gk^{l_{\gk}} \supset 0$ is a filtration of $\gk$, where $l_{\gk} = n-1$ for $\gk = W(n)$, $l_{\gk} = n-2$ for $\gk = S(n)$, $\gk = S'(n)$, and $\gk = \widetilde{H}(n)$, and  $l_{\gk} = n-3$ for $\gk = {H}(n)$. We also have that $\gk = \oplus_{i=-1}^{l_{\gk}} \gk_i$ is a grading for $\gk = W(n), S(n), \widetilde{H}(n)$ and $H(n)$, and that the graded superalgebra $\Gr(S'(n))$ of $S'(n)$ is isomorphic to $S(n)$. The degree-zero component of $\gk$ is described by the isomorphisms  $W(n)_0 \simeq \gg \gl_n$, $S(n)_0 \simeq S'(n)_0 \simeq \gs \gl_n$, $\widetilde{H}(n)_0 \simeq H(n)_0 \simeq \gs \go_n$. 

There is a natural homomorphism $\Ad: \mbox{\bf Aut}\, \Lambda(n) \to \mbox{\bf Aut}\, W(n) $ given by \\ $\Ad_R(\varphi_R)(D):=  \varphi_R D \varphi_R^{-1}$ for $D \in W(n)(R)$ and $\varphi_R \in \mbox{\bf Aut}\, \Lambda(n) (R)$. Note that in this case $\Ad$ is injective. Indeed, let $\varphi_R \xi_i = \sum_{j}c_{ij} \xi_j + f_i $ where $f_i \in \oplus_{i \geq 2}\Lambda(n)_i(R)$ (there is no constant term in the expression of $\varphi_R \xi_i$ because $\xi_i^2 = 0$). Then applying the identity $\varphi_R D f = D \varphi_R f$ for $D = \frac{\partial}{\partial \xi_j}$ and $f = \xi_i$ we find $\frac{\partial}{\partial \xi_j} f_i = 0$ and thus $f_i = 0$. Applying the same identity for $D \in W(n)_0(R) \simeq \gg \gl_n(R)$ we conclude by Schur's Lemma (Lemma \ref{schur}) that $(\varphi_R)_{| \Lambda(n)_1} = \lambda \Id$ for some $\lambda \in R^{\times}$, and finally we let again $D \in W(n)_{-1}(R)$ to find $\lambda = 1$. Since any element $\varphi$ of $\mbox{\bf Aut}\, \Lambda(n) $ is completely determined by its restriction  $\varphi_{| \Lambda(n)_1}$ the kernel of $\Ad$ is trivial as desired. 

It is clear that each element $\theta_R$ of $\mbox{\bf GL}_n(R) = \mbox{\bf Aut}_{R-mod}(R\xi_1 \oplus...\oplus R \xi_n)$ extends to an automorphism of the $R$-algebra $\Lambda(n)\otimes_k R$ (by the universal nature of $\Lambda$). This yields an embedding $\mbox{\bf GL}_n \subset \mbox{\bf Aut}\, \Lambda(n)$. Let $\mbox{\bf N}_{W(n)}$ be the subgroup of $\mbox{\bf Aut}\, \Lambda(n)$ defined as follows:
$$
\mbox{\bf N}_{W(n)}(R):= \{ \nu \in \mbox{\bf Aut} \Lambda(n)(R)\; | \; (\nu- \Id)(\Lambda(n)_i(R)) \subset \oplus_{j>i}\Lambda(n)_j(R), \mbox{ for all }0 \leq i \leq n  \}.
$$
It is clear that $\mbox{\bf N}_{W(n)}$ is a closed (hence affine) subgroup of $\mbox{\bf Aut} \Lambda(n)$. Set for simplicity $\mbox{\bf N}:= \mbox{\bf N}_{W(n)}$, $V := \Lambda(n)_1 = \oplus_{i=1}^n k\xi_i$ and $U := \oplus_{i \geq 1} \Lambda(n)_{2i+1} = \oplus_{i\geq 1}\Lambda^{2i+1}V$.

\begin{lemma} \label{nil-w}
$\mbox{\bf N} \simeq \mbox{\bf Hom}(V,U)$ as $k$-functors.
\end{lemma}
\noindent
{\bf Proof.} If $\nu \in \mbox{\bf N}(R)$ then clearly $\nu - \Id: V(R) \to U(R)$ (remember that $\nu$ preserves $\Lambda(n)_{\bar{1}}$). Since $\nu$ is completely determined by its values on $V(k) = V$ the above yields an injective natural transformation $ \mbox{\bf N} \to \mbox{\bf Hom}(V,U)$. Assume now that an $R$-linear map $\varphi: V(R) \to U(R)$ is given. We show that $\varphi$ is in the image of $\mbox{\bf N}$.  Let $f_i^{(1)}:= \varphi(\xi_i) \in U(R)$. Then $\delta_1:= \sum_{i= 1}^n f_i^{(1)}\frac{\partial}{\partial \xi_i}$  is an $n$-step nilpotent derivation of $\Lambda(n)(R)$. Thus $e^{\delta_1}:= \sum_{i=0}^n \frac{\delta_1^i}{i!} \in  \mbox{\bf Aut} \Lambda(n)(R)$ and $e^{\delta_1}(\xi_i) = \xi_i + f_i^{(1)} +  f_i^{(2)}$ for some $ f_i^{(2)} \in \oplus_{i\geq 2}\Lambda^{2i+1}V(R)$. We next proceed by induction setting  $e^{\delta_{j-1}}...e^{\delta_{1}}(\xi_i) = \xi_i +  f_i^{(1)} +  f_i^{(j)}$ for some $ f_i^{(j)} \in \oplus_{i\geq j}\Lambda^{2i+1}V(R)$ and $\delta_j := - \sum_{i=1}^n f_i^{(j)} \frac{\partial}{\partial \xi_i}$. It is clear that there exists $l$ for which $f_i^{(l)} \neq 0$ for some $i$ and $ f_i^{(j)} = 0$ for all $j >l$ and all $i$. Then $e^{\delta_{l}}...e^{\delta_{1}} - \Id  = \varphi$, and the claim follows. \hfill $\square$

\medskip
It follows from this that as an affine scheme $\mbox{\bf N} = \Spec k[\mbox{\bf N}]$ where $k[\mbox{\bf N}] \simeq S_k(V^* \otimes U)$ ($S_k$ stands for the symmetric algebra over $k$). In particular $\mbox{\bf N}$ is connected. It is also clear by looking at the identical representation of  $\mbox{\bf N}$ on $\Lambda(n)$ that $\mbox{\bf N}$ is unipotent. Furthermore, with the language of \cite{DG}, Ch. IV, \S 2, 4.5, we see that  $\mbox{\bf N} $ is the unipotent group that corresponds to the nilpotent subalgebra $(W(n)^2)_{\bar{0}}$ of  $W(n)$. (In fact one can use this approach to show that  $\mbox{\bf N}$ is connected).

Using the definitions of the Cartan type series $S$, $S'$, and $H$ we easily verify that the following homomorphisms defined by $\Ad$ take place:  
\begin{equation} \label{lin}
\mbox{\bf GL}_n \to  \mbox{\bf Aut}\, S(n), \, \mbox{\bf SL}_{2n} \to \mbox{\bf Aut}\, S'(2n),\,  \mbox{\bf O}_n \times \mbox{\bf G}_m \to \mbox{\bf Aut}\, H(n) 
\end{equation}
 where $ \mbox{\bf G}_m$ corresponds to the scalar matrices $\lambda I_n$, $\lambda \in \mbox{\bf G}_m(R) \simeq R^{\times}$.  All three homomorphisms come from linear maps $\Lambda(n)_1 \to \Lambda(n)_1$ and thus preserve the gradings of $\Lambda(n)$ and $W(n)$. The kernel of the last homomorphism is isomorphic to $\boldsymbol{\mu}_2$ coming from  $( - I_n, - I_n )$ and all other kernels are trivial.

Define the subgroups $\mbox{\bf N}_{S(n)}$,  $\mbox{\bf N}_{S'(n)}$, and $\mbox{\bf N}_{H(n)}$ of $\mbox{\bf Aut}\, \Lambda(n)$ by the identities $\mbox{\bf N}_{S(n)}(R) := \{ \nu \in \mbox{\bf N}(R) \; | \; \widetilde{\nu}^{\theta}(\Delta_n) = \Delta_n\} $, $\mbox{\bf N}_{S'(n)}(R) := \{ \nu \in \mbox{\bf N}(R) \; | \; \widetilde{\nu}^{\theta}(\Delta_n') = \Delta_n'\} $, $\mbox{\bf N}_{H(n)}(R) := \{ \nu \in \mbox{\bf N}(R) \; | \; \widetilde{\nu}^{d}(\omega_n) = \omega_n \} $. The following analog of Lemma \ref{nil-w} describes the nature of these three subgroups.

\begin{lemma}
The following isomorphisms of $k$-functors take place:

(i) $\mbox{\bf N}_{S(n)} \simeq \mbox{\bf N}_{S'(n)}\simeq (\oplus_{i\geq 1}(V^* \otimes \Lambda^{2i+1}V)/\Lambda^{2i} V)_a$

(ii) $\mbox{\bf N}_{H(n)} \simeq (\oplus_{i \geq 2} \Lambda^{2i} V)_a$
\end{lemma}
\noindent
{\bf Proof.} The isomorphisms follow from Propositions 3.3.1 and 3.3.6 in \cite{K2} but for the sake of clearness we sketch the proof. There is a natural isomorphism ${\cal D}: \mbox{\bf Hom} (V(R), \Lambda^{2i+1}V(R)) \to W(n)_{2i}(R)$ defined by ${\cal D}: \psi \mapsto  \sum_{i=1}^n \psi(\xi_j) \frac{\partial}{\partial \xi_i}$.  
With some abuse of notation we will write simply $\widetilde{\psi}^{\theta}$ and  $\widetilde{\psi}^{d}$ instead of $\widetilde{{\cal D}(\psi)}^{\theta}$ and  $\widetilde{{\cal D}(\psi)}^{d}$.

$(i)$ Let $\varphi \in \mbox{\bf Hom} (V(R), \Lambda^{2i+1}V(R))$ be defined by $\varphi(\xi_j) = f_j \xi_j + g_j$, where $f_j, g_j \in \Lambda(n)(R)$ and no summand $c_{j_1,...,j_{2i+1}}\xi_{j_1}...\xi_{j_{2i+1}}$ of $g_j = \sum_{(j_1,...,j_{2i+1})} c_{j_1,...,j_{2i+1}}\xi_{j_1}...\xi_{j_{2i+1}}$ contains $\xi_j$. Then we present $\varphi$ in the following way: $\varphi = \alpha + \beta$, where $\alpha(\xi_j) = \frac{1}{n}(\sum_l f_l)\xi_j $ and $\beta(\xi_j) = g_j + (f_j - \frac{1}{n}(\sum_l f_l))\xi_j $. Note that $\widetilde{\alpha}^{\theta} (\Delta_n) =  (\sum_l f_l) \Delta_n$ and  $\widetilde{\beta}^{\theta} \Delta_n =  0$ (similar identities for $\Delta_n'$). This presentation of $\varphi$ leads to the following isomorphism: 
$$
\mbox{\bf Hom} (V(R), \Lambda^{2i+1}V(R)) \simeq \Lambda^{2i} V(R) \oplus {\cal D}^{-1}(S(n)_{2i}(R)) 
$$
where the first summand corresponds to the subset of  $\mbox{\bf Hom} (V(R), \Lambda^{2i+1}V(R))$ spanned by the inner derivations $\mbox{i}_f: V \to \Lambda^{2i+1} V, \, \mbox{i}_f(g):=fg, f \in \Lambda^{2i} V$. Now using arguments similar to those in the proof of Lemma \ref{nil-w} we show $(i)$.

$(ii)$ The homomorphism $\Lambda^{2i+2} V(R) \to \mbox{\bf Hom} (V(R), \Lambda^{2i+1} V(R))$ defined  by $f \mapsto D_f$ is injective for $i \geq 1$. Its image is $\widetilde{H}(n)_{2i}$ (by the definition of $H$) and we obtain the isomorphism $\Lambda^{2i+2}V(R) \simeq  \{ \nu \in \mbox{\bf Hom} (V(R), \Lambda^{2i+1} V(R)) \; | \; \widetilde{\nu}^d\omega_n = 0\}$ from which the statement easily follows. \hfill $\square$

\medskip
{\bf Note.} Let $\gg$ be one of the Lie superalgebras $W(n), S(n)$, $S'(n)$ (if $n=2l$), $H(n)$, or $\widetilde{H}(n)$. Denote by $\gn_{\gg}$ the nilpotent Lie algebra $(\gg^2)_{\bar{0}}$. Note that $\gn_{\gg}$ is the radical of the Lie algebra $\gg_{\bar{0}}$ for $\gg = S(n), S'(2l), H(n), \widetilde{H}(n)$ and it is the commutator of the radical of  $\gg_{\bar{0}}$ for $\gg = W(n)$.  We can easily show that $\mbox{\bf N}_{\gg}$ is the unipotent group that corresponds to $\gn_{\widetilde{\gg}}$, where $\widetilde{\gg} := \widetilde{H}(2l)$ for $\gg = H(2l)$ and $\widetilde{\gg} := \gg$ in all other cases (the arguments for $\gg = S, S', H$ are similar to those for $\gg = W$).
\medskip

From the discussion above we see that $\Ad : \mbox{\bf N}_{\gg} \to \mbox{ \bf Aut}\,\gg$ is an injective homomorphism. Combining this homomorphism with the homomorphisms in (\ref{lin}) we find the following injective homomorphisms:

$\mbox{\bf N}_{W(n)} \rtimes \mbox{\bf GL}_n \to \mbox{ \bf Aut}\, W(n), \, \mbox{\bf N}_{S(n)} \rtimes \mbox{\bf GL}_n \to \mbox{ \bf Aut}\, S(n)$, 

$\mbox{\bf N}_{S'(2l)} \rtimes \mbox{\bf SL}_{2l}  \to \mbox{ \bf Aut}\, S'(2l), \, \mbox{\bf N}_{H(n)} \rtimes ((\mbox{\bf O}_{n} \times  \mbox{ \bf G}_m) /\boldsymbol{\mu}_2) \to \mbox{ \bf Aut}\, H(n).$

We easily see that $(\mbox{\bf O}_{n} \times  \mbox{ \bf G}_m) /\boldsymbol{\mu}_2 \simeq \mbox{\bf SO}_{n} \times  \mbox{ \bf G}_m$ for odd $n$. In the case when $n$ is even one reasons as in \S \ref{osp} and shows that there is a split exact sequence:
$$
1 \to \mbox{\bf N}_{H(n)} \rtimes ((\mbox{\bf SO}_{n} \times  \mbox{ \bf G}_m) /\boldsymbol{\mu}_2) \to \mbox{\bf Aut}\, H(n) \to \Z_{2,k} \to 1
$$
where $\Z_{2,k} = <r_m>_k$ for $r_m \in \mbox{\bf O}_{n}$, $\det r_m = -1$, and $r_m^2 = 1$. As usual, using Serganova's results that the corresponding homomorphisms  are surjective on $k$-points we show that all four homomorphisms above are isomorphisms.

\section{The supercentroid of $\mathfrak g$} \label{centro}
\setcounter{equation}{0}

In this section we find the {\it supercentroid} $\Ctd(\gg)$ of $\gg$, where 
$$
\Ctd_k(\gg):=\{\chi \in \End_{k}(\gg) \; | \; \chi([x,y]) = [\chi(x),y] \mbox{ for all } x,y \in \gg\}.
$$
Recall that by definition of $\End_{k}(\gg)$, we have $\varphi(\gg_{\bar{i}}) \subseteq \gg_{\bar{i}}$ for all $\varphi \in \End_{k}(\gg)$ and $\bar{i} \in \Z_2$. Let $\lambda_{\gg}: k \to \End_{k}(\gg)$ be defined by $\lambda_{\gg}(a)(x):=ax$ for $a \in k$ and $x \in \gg$. We call $\gg$ {\it central} if $\Ctd_k(\gg) = \lambda_{\gg}(k)$.

\begin{proposition} \label{central}
All simple finite dimensional Lie superalgebras are central.
\end{proposition}
\noindent
{\bf Proof.} It is clear that an endomorphism $\chi$ is in the centroid of $\gg$ if and only if $\chi \circ \ad = \ad \circ \chi$. The proposition follows directly from the superversion of the classical Schur's Lemma (see \S 1.1.6 in \cite{K2}) since the adjoint representation of $\gg$ is irreducible. \hfill $\square$

\section{The Proofs of Theorems \ref{cohomo}, \ref{th3}, \ref{th3.5}, and \ref{th4}} \label{proof2}

\setcounter{equation}{0}
\subsection{Cohomological considerations} \label{cohocon}

We refer the reader to \cite{Mil}, \cite{DG}, and [SGA] for general facts about \'etale nonabelian $H^1$. Let $L(\gg, \sigma)$ be as in \S \ref{loop} with $\sigma$ of period $m$. Recall that for $S = k[z^{\pm 1}]$ viewed as an $R = k[t^{\pm 1}]$-algebra via $t \mapsto z^m$ we have
$$
L(\gg, \sigma) \otimes_R S \simeq_S \gg \otimes S \simeq (\gg \otimes_k R)\otimes_R S.
$$
Because of this, we can attach to $L(\gg, \sigma)$ an element 
$$
[L(\gg, \sigma)] \in H_{\acute et}^1(S/R, \mbox{\bf Aut}\, \gg_R) \subset H_{\acute et}^1(R, \mbox{\bf Aut}\, \gg_R).
$$
Now $S/R$ is a Galois extension with Galois group $\Gamma \simeq \Z_m$ given by $\Gamma = <\gamma>$ where $\gamma(z) = \zeta_m z$. One knows that there is a natural correspondence (see \cite{W} or \cite{Mil}, \S III.2, Example 2.6) $ H_{\acute et}^1(S/R, \mbox{\bf Aut}\, \gg_R) \stackrel{\sim}\to H^1(\Gamma, \mbox{\bf Aut}\, \gg_R(S))$ where this last is the ``usual'' nonabelian cohomology. Here $\Gamma$ acts on $\mbox{\bf Aut}\, \gg_R(S):= \Aut_S (\gg \otimes R) \otimes_S S = \Aut_S \gg \otimes S$ by ``acting on each of the coordinates'' (the elements of $\mbox{\bf Aut}\, \gg_R(S)$ can, after fixing a homogeneous basis of $\gg$, be thought as $n \times n$ matrices with entries in $S$ where $n = \dim_k \gg$). In terms of maps, under this action  
${}^\gamma \theta : \gg \otimes S \to \gg \otimes S$  is given by ${}^\gamma \theta := (1 \otimes \gamma) \circ \theta \circ (1\otimes -\gamma)$.
 The action of $\Gamma$ can of course also be described functorially as follows: If $\theta \in \mbox{\bf Aut}\, \gg_R(S) = \Hom_{R\mbox{-}alg}(R[\mbox{\bf Aut}\, \gg_R], S)$, then ${}^{\gamma}\theta = \gamma \circ \theta$.
An obvious calculation shows that the unique homomorphism $u: \Gamma \to \mbox{\bf Aut}\, \gg_R(S):= \mbox{Aut}_S(\gg \otimes_k S)$ given by $\gamma \mapsto \sigma^{-1} \otimes 1$ is in fact a cocycle in $Z^1(\Gamma, \mbox{\bf Aut}\, \gg_R(S))$. Then $[L(\gg, \sigma)]$ is nothing but the cohomology class of this cocycle.
 
Via our choice of compatible roots of unity, the algebraic fundamental group of $R$ based at the geometric point $1$ can be identified with $\varprojlim \Z_n$, hence cyclic in the topological sense. From Grothendieck's work we obtain the correspondences 
$$
H_{\acute et}^1(R, \mbox{F}_R ) \simeq \mbox{Conjugacy classes of F}  \simeq \mbox{Galois extensions of }R\mbox{ with group F}.
$$
(see \cite{SGA1}, \S 5, Exp. 11 and \cite{P}). As it seems inevitable, under the canonical map $H_{\acute et}^1(R, \mbox{\bf Aut}\, \gg_R ) \to H_{\acute et}^1(R, \mbox{F}_R )$ we have $[L(\gg, \sigma)] \mapsto \{\mbox{conjugacy class of } \bar{\sigma}^{-1}\}$ where 
$\, \bar{} : \mbox{\bf Aut}\, \gg (k)\to \mbox{F}$ arises from Theorem \ref{th1}.

\subsection{Proof of Theorem \ref{cohomo}} \label{proof2.2}

That the map $ H_{\acute et}^1(R, \mbox{\bf Aut}\, \gg_R ) \to H_{\acute et}^1(R, \mbox{F}_R )$ is surjective follows from the fact that in all cases (as made clear in the proof of Theorem \ref{th1}) the group $\mbox{F}_k$ lifts to a finite constant subgroup  $\widetilde{\mbox{F}}_k$ of $\mbox{\bf Aut}\, \gg$, and one can now conclude by the considerations of \S \ref{cohocon}.

Let $\pi \in \mbox{F}$, and lift $\pi$ to an element $\widetilde{\pi}$ of $\widetilde{\mbox{F}}$. Then $\widetilde{\pi}$ induces an automorphism  $\widetilde{\pi}_0$ of  $\mbox{\bf G}^0$ and hence by base change also one of  $\mbox{\bf G}^0_R$. One can then consider the twisted group   ${}^{\widetilde{\pi}}\mbox{\bf G}^0_R$ constructed as follows. Let $l$ be the order of $\widetilde{\pi}_0$. Consider $S = k[z^{\pm 1}]$ which we view as an $R$-algebra via $t \mapsto z^l$. As we have seen $S/R$ is a Galois extension with Galois group $\Gamma \simeq \Z_l$ generated by $\gamma$ where $\gamma(z) = \zeta z$ where $\zeta = \zeta_l$.

If $A \in R${\it -alg} and $x \in \mbox{\bf G}_R^0(S \otimes_R A):= \Hom_R(R[\mbox{\bf G}^0], S\otimes_R A)$ we define ${}^\gamma x \in \mbox{\bf G}_R^0(S \otimes_R A)$ by ${}^\gamma x = (\gamma \otimes 1) \circ x$. Then the twisted group ${}^{\widetilde{\pi}}  \mbox{\bf G}_R^0$ is the $R$-group whose functor of points is given by ${}^{\widetilde{\pi}} \mbox{\bf G}_R^0 (A) = \{ x \in \mbox{\bf G}_R^0 (S \otimes_R A)\; : \; \widetilde{\pi}_{0_{S \otimes_R A}}({}^\gamma x) = x\}$.

\medskip
\noindent
{\bf Remark.} That ${}^{\widetilde{\pi}} \mbox{\bf G}_R^0$ is a an $R$-group (i.e. affine) follows in general from descent. In the present case one can easily see what the coordinate ring is by viewing $\pi$ as an automorphism of the Hopf algebra $R[\mbox{\bf G}^0] = R \otimes_k k[\mbox{\bf G}^0]$. Then $R[{}^{\widetilde{\pi}} \mbox{\bf G}_R^0] = \{f \in S \otimes_R R[\mbox{\bf G}]\; | \; \widetilde{\pi}_{0_S}(\gamma \otimes 1)(f) = f \}$. (This amounts to the usual cocycle condition of descent because $\Gamma$ is generated by $\gamma$ and the action of $\Gamma$ commutes with that of $\widetilde{\pi}$).  Up to an isomorphism the twisted group  ${}^{\widetilde{\pi}} \mbox{\bf G}_R^0$ does not depend on the choice of the lift $\widetilde{\pi}$ and will henceforth simply be denoted by ${}^{\pi} \mbox{\bf G}_R^0$. 

\smallskip
We now come to the most delicate point of the proof, namely showing that the map $H_{\acute et}^1(R, \mbox{\bf Aut}\, \gg_R ) \to  H_{\acute et}^1(R, \mbox{F}_R )$ is injective. The complication stems from the fact that these $H^1$'s are not groups, but rather sets with a distinguished element. As a consequence, to show injectivity it is not enough that $H_{\acute et}^1(R, \mbox{\bf G}^0_R ) = 0$, but rather that {\it all} twisted $ H_{\acute et}^1(R, {^\pi}\mbox{\bf G}_R^0 )$ also vanish. Unlike $ \mbox{\bf G}^0_R$, which is obtained from the $k$-group  $ \mbox{\bf G}^0$ by base change, the twisted  ${}^\pi \mbox{\bf G}_R^0$  (unless $\pi = 1$) {\it do not come from any algebraic} $k$-group (there are no twisted forms of  $ \mbox{\bf G}^0$ over $k$). Thus   ${}^\pi \mbox{\bf G}_R^0$ is a reductive group over $R$ in the sense of \cite{SGA3}. Computing its cohomology is therefore not surprisingly much more delicate than in the untwisted case.

The general strategy to establish Theorem \ref{cohomo} will be as follows. We will consider an exact sequence of the form
\begin{equation}\label{innauto}
1 \to \mbox{\bf z} \to  \mbox{\bf G} \to  \mbox{\bf G}^0 \to 1.
\end{equation}
where $\mbox{\bf G} $ is reductive and $\mbox{\bf z}$ is central and finite, i.e. an isogeny.  We also assume that this sequence can be twisted by $\pi$ to yield
$$
1 \to {}^\pi\mbox{\bf z}_R \to  {}^\pi\mbox{\bf G}_R \to  {}^\pi\mbox{\bf G}^0_R \to 1.
$$
Passing to cohomology we obtain
$$
\to H_{\acute et}^1(R, {}^\pi\mbox{\bf G}_R) \to  H_{\acute et}^1(R, {}^\pi\mbox{\bf G}^0_R) \to  H_{\acute et}^2(R, {}^\pi\mbox{\bf z}_R).
$$
By \cite{SGA4}, Exp. IX, 5.7, and Exp. X, 5.2 our ring $R$ is of cohomological dimension $1$. As a consequence $H_{\acute et}^2(R, {}^\pi\mbox{\bf z}_R) = 1$ and therefore Theorem \ref{cohomo} will follow once we establish that in each case for all $\pi \in \mbox{F}$ 
\begin{equation}\label{cohovan}
H_{\acute et}^1(R, {}^\pi\mbox{\bf G}_R) = 1.
\end{equation}

\medskip
\noindent
{\bf \ref{proof2.2}.1 The case} $\mbox{F} = \mbox{1}$. Let $\mbox{\bf z} = {1}_k$ in (\ref{innauto}). Then $\mbox{\bf Aut}\, \gg = \mbox{\bf G}^0$ is a connected linear algebraic $k$-group. By \cite{P}, Proposition 5, $H_{\acute et}^1(R,\mbox{\bf G}_R^0 ) = \mbox{1}$.

In view of this we may therefore assume that $\mbox{F} \neq \mbox{1}$, that $\pi \neq 1$, and verify (\ref{cohovan}) in each of the cases.

\medskip
\noindent
{\bf \ref{proof2.2}.2 The case $\gg = \gs \gl (m|n)$}. The generator of the outer $\Z_2$ lifts to the element $r$ of $\mbox{\bf Aut}\, \gs \gl (m|n)(k)$, which acts on each of the factors of $\mbox{\bf SL}_m \times \mbox{\bf SL}_n \times \mbox{\bf G}_m $ via $\pi_R : X \mapsto (X^{t})^{-1}$ in the $\mbox{\bf SL}$'s and $\gamma_R : X \mapsto X^{-1}$ in  $\mbox{\bf G}_m$. This action stabilizes $\boldsymbol{\mu}_m \times \boldsymbol{\mu}_n$ where it acts via $(\alpha, \beta) \mapsto (\alpha^{-1}, \beta^{-1})$. The exact sequence  (\ref{innauto}) is given by
$$
1 \to \boldsymbol{\mu}_m \times \boldsymbol{\mu}_n \to \mbox{\bf SL}_m \times \mbox{\bf SL}_n \times \mbox{\bf G}_m \to \mbox{\bf SL}_m \times \mbox{\bf SL}_n \times \mbox{\bf G}_m /(\boldsymbol{\mu}_m \times \boldsymbol{\mu}_n) \to 1.
$$
To establish (\ref{cohovan}) we need 
\begin{equation}\label{coho-sl}
H_{\acute et}^1(R, {}^{\gamma} \mbox{\bf G}_{m_R} ) = H_{\acute et}^1(R, {}^{ \pi} \mbox{\bf SL}_{l_R} )  = 1
\end{equation}

Let $S = k[z^{\pm 1}]$ which we view as an $R$-algebra via $t \mapsto z^2$. This is a Galois extension of degree $2$. The norm map $N: S^{\times} \to R^{\times}$ induces a surjective $R$-group homomorphism
$$
N: {\cal R}_{S/R}  \mbox{\bf G}_{m_S} \to  \mbox{\bf G}_{m_R}
$$
where $ {\cal R}$ denotes Weyl restriction. A direct calculation shows that our twisted multiplicative group is the kernel of $N$. Thus
$$
1 \to {}^{\gamma}\mbox{\bf G}_{m_R} \to  {\cal R}_{S/R}  \mbox{\bf G}_{m_S} \to \mbox{\bf G}_{m_R} \to 1.
$$
By passing to cohomology we get 
$$
S^{\times} \stackrel{N}{\to} R^{\times} \to H_{\acute et}^1(R, {}^{\gamma}\mbox{\bf G}_{m_R}) \to H_{\acute et}^1(R,{\cal R}_{S/R}  \mbox{\bf G}_{m_S} ).
$$
The last $H^1$ vanishes by Shapiro's Lemma since 
$$
 H_{\acute et}^1(R,{\cal R}_{S/R}  \mbox{\bf G}_{m_S} ) = H_{\acute et}^1(S, {}^{\gamma}\mbox{\bf G}_{m_S}) = \mbox{Pic}(S) = 1.
$$
The norm map is given by $\lambda z^l \mapsto (-1)^l \lambda^2t^l$, hence it is surjective since $k$ is algebraically closed. It follows that $H_{\acute et}^1(R, {}^{\gamma}\mbox{\bf G}_{m_R} ) = 1$ as desired.

Finally ${}^{\pi}\mbox{\bf SL}_{l_R}$ is quasisplit of simply connected type so its $H^1$ is trivial by \cite{P}, Proposition 3. This finishes the proof of (\ref{coho-sl}).

\medskip
\noindent
{\bf \ref{proof2.2}.3 The case $\gg = \gp \gs \gl (n|n), n>2$}. The reasoning is similar to that of $\gs \gl (m|n)$ except that we now have more twisted groups to content with. If $\sigma$ denotes the element of $\mbox{\bf Aut}\, (\mbox{\bf SL}_n \times \mbox{\bf SL}_n)$ that comes from ``switching'' the two factors, we must show in addition to (\ref{coho-sl}) that
\begin{equation}\label{coho-psl}
\begin{array}{cccccl} 
H_{\acute et}^1(R, {}^\sigma (\mbox{\bf SL}_n \times \mbox{\bf SL}_n )_R) & = & H_{\acute et}^1(R, {}^{\sigma \pi} (\mbox{\bf SL}_n \times \mbox{\bf SL}_n )_R) & = & 1. &  \\
\end{array}
\end{equation}
Again both $ {}^{\sigma} (\mbox{\bf SL}_n \times \mbox{\bf SL}_n )_R$ and  $ {}^{\sigma \pi} (\mbox{\bf SL}_n \times \mbox{\bf SL}_n )_R$ are quasisplit of simply connected type so their $H_{\acute et}^1(R, -)$ are trivial.

\medskip
\noindent
{\bf \ref{proof2.2}.4 The case  $\gg = \gp \gs \gl (2|2)$.} Take $\mbox{\bf z} = \mbox{1}_k$ in (\ref{innauto}). The cohomology that needs to vanish is
$H_{\acute et}^1(R,{}^\gamma \mbox{\bf G}_m)$,  
and we have already seen that this is indeed the case.

\medskip
\noindent
{\bf \ref{proof2.2}.5 The case $\gg = \gp \gs \gq (n)$.} Take $\mbox{\bf z} = \mbox{1}_k$ in (\ref{innauto}). All twisted forms of $\mbox{\bf PGL}_{n_R}$ are $R$-groups of adjoint type and hence have trivial $H_{\acute et}^1(R, -)$ by Proposition 3 of \cite{P}.

\medskip
\noindent
{\bf \ref{proof2.2}.6 The case  $\gg = \go \gs \gp (m|2n)$, $m$ even.} Again we take $\mbox{\bf z} = \mbox{1}_k$ in  (\ref{innauto}). We need 
$$
H_{\acute et}^1(R, {}^r (\mbox{\bf SO}_m \times \mbox{\bf Sp}_{2n} )_R) = 1.
$$
Now, ${}^r (\mbox{\bf SO}_m \times \mbox{\bf Sp}_{2n} )_R \simeq {}^\pi \mbox{\bf SO}_{m_R} \times \mbox{\bf Sp}_{2n_R} $ where $\pi$ is the automorphism of the ``tail'' of the Dynkin diagram of the group $\mbox{\bf SO}_{m}$. One knows that $H_{\acute et}^1(R, {}^\pi \mbox{\bf SO}_m ) = 0$. To see this consider the simply connected cover
$$
1 \to \boldsymbol{\mu}_2 \to \mbox{\bf Spin}_m \to \mbox{\bf SO}_m \to 1
$$
and conclude with the aid of Proposition 3 of \cite{P} by passing to cohomology.

\medskip
\noindent
{\bf \ref{proof2.2}.7 The case  $\gg = D(1)$.} For (\ref{innauto}) we take 
$$
1 \to \boldsymbol{\mu}_2 \times \boldsymbol{\mu}_2 \to \mbox{\bf SL}_2 \times  \mbox{\bf SL}_2 \times  \mbox{\bf SL}_2 \to  \mbox{\bf Aut}\, D(1) \to 1,
$$
which we need to twist by $\pi$ where $\pi$ switches the first two factors of $(\mbox{\bf SL}_2)^3$ and the two factors of $\boldsymbol{\mu}_2 \times \boldsymbol{\mu}_2$. The twisted $(\mbox{\bf SL}_2)_R^3$ is of simply connected type and hence has trivial $H_{\acute et}^1(R, -)$.

\medskip
\noindent
{\bf \ref{proof2.2}.8 The case  $\gg = D(\alpha), \alpha \neq 1, \alpha^3 = 1$.} The reasoning is similar to that of $D(1)$.

\medskip
\noindent
{\bf \ref{proof2.2}.9 The case $\gg = H(2l)$.} The reason is similar to that of $\gg = \go \gs \gp (2l|2n)$. To deal with the twisted unipotent group, one reasons by d\'evissage on a suitable central composition series. 

\subsection{Proof of Theorem \ref{th3}}
The plan is to describe how the results we have obtained thus far allow us to apply the central ideas of \cite{ABP} to obtain Theorem \ref{th3}

First we observe that because $\gg$ is perfect and central (see Proposition \ref{central}), $\Ctd_k(L(\gg, \sigma)) \simeq \Ctd_R(L(\gg, \sigma)) \simeq R$ (reason as in \cite{ABP}, Lemma 4.3). Now if $\varphi: L(\gg_1, \sigma_1) \to L(\gg_2, \sigma_2)$ is an isomorphism of $k$-Lie superalgebras, there exists a unique $\widetilde{\varphi} \in \mbox{\bf Aut}_k(R)$ such that $\varphi(rx) = \widetilde{\varphi}(r)\varphi(x)$ (ibid Lemma 4.4). To say that $\varphi$ is an isomorphism of $R$-Lie superalgebras is equivalent to saying that $\widetilde{\varphi} = \Id$. 

The trick now is to try to replace $\varphi$ by another isomorphism $\varphi'$ which is $R$-linear (there may be a price for doing this as we will see). We have $\mbox{Aut}_k(R) \simeq k^{\times} \rtimes \Z_2$ where $\lambda \in  k^{\times}$ accounts for $t \mapsto \lambda t$ and the generator of $\Z_2$ for the switch $t \to t^{-1}$. The up shot of this is that 

\begin{equation} \label{sigmas}
L(\gg_1, \sigma_1) \simeq_k L(\gg_2, \sigma_2) \Leftrightarrow L(\gg_1, \sigma_1) \simeq_R L(\gg_2, \sigma_2) \mbox{ or } L(\gg_1, \sigma_1) \simeq_R L(\gg_2, \sigma_2^{-1}).
\end{equation}
(ibid Theorem 4.6).

The reason behind this artifice is to bring us to the situation where the loop algebras are viewed as $R$-Lie superalgebras, and as such, their isomorphism classes are completely understood by cohomological methods. Assume $L(\gg_1, \sigma_1) \simeq_R L(\gg_2, \sigma_2)$ with $\sigma_1$ and $\sigma_2$ of common period $m$ (the case of $\sigma_2^{-1}$ is similar). Then as $S$-Lie superalgebras (see \S \ref{cohocon})
$$
\gg_1 \otimes_k S \simeq L(\gg_1, \sigma_1)\otimes_R S \simeq L(\gg_2, \sigma_2)\otimes_R S \simeq \gg_2 \otimes_k S.
$$
Now tensoring over $S$ with the algebraic closure $K$ of the quotient field  of $S$ we obtain $\gg_1 \otimes_k K \simeq \gg_2 \otimes_k K $ as $
k$-Lie superalgebras and hence $\gg_1 \simeq \gg_2$ by the classification of simple Lie superalgebras over algebraically closed fields of characteristic zero. \hfill $\square$

\subsection{Proof of Theorem \ref{th3.5}}

Assume $L(\gg, \sigma_1) \simeq_R L(\gg, \sigma_2)$. Then $[L(\gg, \sigma_1)] = [L(\gg, \sigma_2)]$ and hence by \S \ref{cohocon} we obtain $\sigma_1^{-1} \sim \sigma_2^{-1}$ hence $\bar{\sigma}_1 \sim \bar{\sigma}_2$ in $\mbox{F}$. Similarly, if $L(\gg, \sigma_1) \simeq_R L(\gg, \sigma_2^{-1})$ then $\bar{\sigma}_1 \sim \bar{\sigma}_2^{-1}$. That $(i)$ implies $(ii)$ now follows from (\ref{sigmas}). For the converse  observe that if  $\bar{\sigma}_1 \sim \bar{\sigma}_2$ then $[L(\gg, \sigma_1)] = [L(\gg, \sigma_2)]$ by \S \ref{cohocon} in view of Theorem \ref{cohomo}. But then $L(\gg, \sigma_1) \simeq_R L(\gg, \sigma_2)$ and a fortiori these two are isomorphic as $k$-Lie superalgebras. Similarly $\bar{\sigma}_1 \sim \bar{\sigma}_2^{-1}$ implies $L(\gg, \sigma_1) \simeq_k L(\gg, \sigma_2^{-1})$. But since  $L(\gg, \sigma) \simeq_k L(\gg, \sigma^{-1})$ always (\cite{ABP}, Lemma 2.4(a)) the proof is now complete.  \hfill $\square$

\subsection{Proof of Theorem \ref{th4}}

The Zen of torsors shows that the isomorphism classes of forms of the $R$-Lie superalgebra $\gg \otimes_k R$ are classified by $H^1_{\acute et}(R, \mbox{\bf Aut}\, \gg_R)$ (we may replace the flat site by the \'etale site of $\Spec R$ since our $R$-group $\mbox{\bf Aut}\, \gg_R$ is smooth). By applying now Theorem \ref{cohomo} one concludes that our forms are parameterized by $H^1_{\acute et}(R, \mbox{F}_R)$. But the considerations of \S 8.1 show that all elements of $H^1_{\acute et}(R, \mbox{F}_R)$ are of the form $[L(\gg, \widetilde{\pi})]$ for some finite order automorphism $\widetilde{\pi}$ of $\gg$.  \hfill $\square$

Department of Mathematical and Statistical Sciences

University of Alberta

Edmonton, Alberta T6G 2G1

CANADA
\medskip

E-mails: grant@math.ualberta.ca, a.pianzola@math.ualberta.ca

\end{document}